\documentclass{article}
\usepackage[latin1]{inputenc} 
\usepackage[T1]{fontenc}      
\usepackage[frenchb,english]{babel} 
\usepackage{amsmath}
\usepackage{amssymb}
\usepackage{amsthm}
\makeatletter
\@addtoreset{equation}{section}         

\newcommand {\N}{\mathbb{Z}_{+}}  
\newcommand {\Z}{\mathbb{Z}}    
\newcommand {\R}{\mathbb{R}}  
\newcommand {\C}{\mathbb{C}}  
\newcommand {\D}{\mathbb{D}}  
\newcommand {\T}{\mathbb{T}}  

\swapnumbers
\theoremstyle{plain}
\newtheorem{theorem}{Theorem}[section]
\newtheorem{prop}[theorem]{Proposition}
\newtheorem{lem}[theorem]{Lemma}
\newtheorem{cor}[theorem]{Corollary}

\theoremstyle{definition}
\newtheorem{definition}[theorem]{Definition}
\newtheorem{remark}[theorem]{Remark}
\newtheorem{example}[theorem]{Example}

\newcommand{\Ker}{{\mbox{Ker }}}

\newcommand{\LL}{\mathcal{B}}
\renewcommand{\H}{\mathcal{H}}

\newcommand{\cd}{\cdots}

\newcommand{\noi}{\noindent \smallskip }
\newcommand{\beq}{\begin{equation}}
\newcommand{\eeq}{\end{equation}}
\newcommand{\beqr}{\begin{eqnarray*}}
\newcommand{\eeqr}{\end{eqnarray*}}
\newcommand{\bc}{\begin{center}}
\newcommand{\ec}{\end{center}}

\newcommand{\gm}{\gamma}

\newcommand{\ep}{\varepsilon}

\newcommand{\om}{\omega}

\def\firsta{\  {\text {\rm (a)}}\ \ }

\def\condb{\smallskip  \noindent {\text {\rm (b)}}\ \ }
\def\condc{\smallskip  \noindent {\text {\rm (c)}}\ \ }

\begin{document}

\title{Constrained von Neumann
inequalities
\thanks{Preprint. ; to appear in \emph{Advances in Mathematics}}}

\author{{\bf C. BADEA} \\
$\; \; $\\
D\'epartement de Math\'ematiques, UMR 8524 au CNRS\\  
Universit\'e de Lille I, F-59655 Villeneuve d'Ascq, France\\
{\tt Catalin.Badea@math.univ-lille1.fr}
$\; \; $\\
$\; \; $\\
and
$\; \; $\\
$\; \; $\\
{\bf G. CASSIER}\\
$\; \; $\\
Institut G. Desargues, UMR 5028 au CNRS\\
Universit\'e de Lyon I, F--69622 Villeurbanne, France\\
{\tt Gilles.Cassier@desargues.univ-lyon1.fr}}
\date{}
\maketitle

\begin{abstract}
An equivalent formulation of the von Neumann inequality 
states that the
backward shift $S^*$ on $\ell_{2}$ is extremal, in the sense that if
$T$ is a Hilbert space contraction, then $\|p(T)\|
\leq \|p(S^*)\|$ for each polynomial $p$.
We discuss several results of the following
type~: if $T$ is a Hilbert space contraction satisfying some
constraints, then $S^*$ restricted to a suitable invariant subspace is
an extremal operator. Several operator radii are used
instead of the operator norm. Applications to inequalities of
coefficients of rational functions positive on the torus are given.
\end{abstract}

\selectlanguage{frenchb}
\begin{abstract}
 D'apr\`es l'in\'egalit\'e de von Neumann, l'op\'erateur $S^*$
de translation en arri\`ere est extr\'emal : on a $\|p(T)\|
\leq \|p(S^*)\|$ pour chaque contraction hilbertienne $T$ et chaque polyn\^ome
$p$. Nous d\'emontrons que, pour les contractions hilbertiennes v\'erifiant
certaines contraintes, la restriction de $S^*$ \`a un sous-espace invariant est
un op\'erateur extr\'emal.
\end{abstract}
\vspace{3mm}
\selectlanguage{english}


\vspace{3mm}

\thispagestyle{empty}
\newpage
\setcounter{section}{-1} 
\section{Introduction}
Let $T$ be a Hilbert space contraction, that is a bounded
linear operator of norm at most one on a complex, separable Hilbert space $H$.
A well-known inequality due to J. von Neumann \cite{vNe} asserts that
\begin{equation}
     \label{vNineq}
    \| p(T) \| \leq \| p \|_{\infty},
\end{equation}
for every polynomial $p\in \C[X]$. Here
$$\| p \|_{\infty} = \sup \{ | p(z) | : z \in \C , | z | \leq 1 \}$$
is the supremum norm of $p$, while
$$\| p(T) \| = \| p(T) \|_{\LL (H)}$$
is the operator norm
of $p(T)$ in $\LL (H)$, the $C^{\ast}$-algebra of all bounded linear
operators on $H$. The same inequality extends for functions in the disc
algebra
$A(\D)$ and, if $T$ is a completely non-unitary (c.n.u.) contraction, it
extends
to bounded analytic functions $f\in H^{\infty}(\D)$ \cite{SNF}. Recall
that a c.n.u. operator is one which has no unitary direct summand \cite{SNF}.

Denote by $S$ the
forward unilateral shift on $\ell^2$,
$$S(x_{0},x_{1}, \cdots ) = (0,x_{0},x_{1},\cdots ) ,$$
and by $S^{\ast} \in B(\ell^2)$,
$$S^{\ast}(x_{0},x_{1}, \cdots ) = (x_{1},x_{2}, \cdots ) ,$$
its adjoint (the backward shift).

An equivalent formulation of the von Neumann inequality (\ref{vNineq})
is the following :
for every Hilbert space contraction $T$ and every polynomial $p$
we have
\begin{equation}
    \label{vNineq2}
    \| p(T) \|_{\LL(H)} \leq \| p(S^{\ast}) \|_{\LL (\ell^2)}.
\end{equation}
We say that $S^*$ is extremal. A proof of the inequality 
(\ref{vNineq2}) will be sketched
in Section 2.

We will discuss several results of the following
type~: if $T$ is a Hilbert space contraction satisfying some
constraints and $\om$ is an operator radius, then there exists a
suitable invariant subspace $E$ of $S^*$ such that
$$\om(p(T)) \leq \om(p(S^*\mid E)),$$
that is $S^*$ restricted to a suitable invariant subspace is
an extremal operator.

Several results of this type are known in the literature.
The following result was proved by V. Pt\'ak
\cite{Pta1}, \cite{Pta2} in a particular case~; the general case
was proved by Pt\'ak and N.J. Young \cite{ptak/young}. Suppose that
$p$ and $q$ are arbitrary analytic polynomials.
    Let $T$ be a Hilbert space contraction of spectral radius smaller than
    one and suppose that $q(T) = 0$. Then
$$\| p(T) \| \leq \| p(S^{\ast}\mid \Ker q(S^{\ast})) \| .$$

The following extension was given by B. Sz.-Nagy \cite{sznagy}.
Let $f$ and $g$ be two functions in $H^{\infty}(\D)$. Let $T$ be a
    Hilbert space c.n.u. contraction such that $g(T) = 0$.
    Then
$$\| f(T) \| \leq \| f(S^{\ast}\mid \Ker g(S^{\ast})) \| .$$
An equivalent form of the Sz.-Nagy's result was stated by J.P.
Williams \cite{williams:MR}~; Williams' proof is given in the survey
paper \cite{ptak:casopis}.

An apparently unrelated inequality due to U. Haagerup and
P. de la Harpe \cite{HaHa}
asserts that each bounded linear
nilpotent
contraction $T$ with $T^n = 0$, $n \geq 2$, satisfies
the inequality
\begin{equation} \label{HHineq}
\om_{2}(T) \leq \cos \frac{\pi}{n+1}.
\end{equation}
Here $\om_{2}(T) $ denotes the numerical radius of $T$ defined by
$$
\om_{2}(T) = \sup \{ | \langle Tx|x\rangle | : x \in H, \| x \| = 1 \}.
$$

To see how the Haagerup-de la Harpe inequality fits into the present
framework,
let $S_n^{*}$ be the nilpotent Jordan cell
$$S_{n}^* = \left(\begin{array}{cccccc}
            0 & 1 & 0 & \cd & 0 & 0 \\
            0 & 0 & 1 & \cd & 0 & 0 \\
            \vdots & \vdots & \vdots & \ddots & \vdots & \vdots \\
            0 & 0 & 0 & \cd & 0 & 1 \\
            0 & 0 & 0 & \cd & 0 & 0
        \end{array}\right) $$
on the standard
Euclidean space $\C^{n}$.
Then \cite{gust/rao:book} $\cos (\pi/(n+1)) = \om_{2}(S_{n}^*)$ and
$S_{n}^*$ is unitarily equivalent to $S^{\ast}\mid \C^n =
S^{\ast}\mid \Ker u_{n}(S^{\ast})$, where $u_{n}(z) = z^n$. Therefore
the inequality of Haagerup and de la Harpe
states that if $u_{n}(T) = 0$, then
$$
\om_{2}(T) \leq \om_{2}(S^{\ast}\mid \Ker u_{n}(S^{\ast})) \; \; .
$$
We refer to \cite{wu,suen,pop} for recent papers related to this inequality.

In \cite{HaHa}, inequality (\ref{HHineq}) is shown to be equivalent to
an inequality, due to L. Fejer (1915), for the first
coefficient $c_{1}$ of a positive trigonometric polynomial
$\sum_{j=-n+1}^{n-1}c_{j}e^{ijt}$, namely
\begin{equation*}
\left| c_{1}\right| \leq c_{0}\cos (\frac{\pi }{n+1}) .
\end{equation*}
We will prove other inequalities for
coefficients of rational functions positive on the torus or for
coefficients of positive trigonometric polynomials which are
related to our constrained von Neumann inequalities. In particular,
we obtain (Theorem \ref{thm:harm}) the following inequality for the
sum of the absolute value of two coefficients of a positive
trigonometric polynomial of degree $n$~:
\begin{equation*}
\left| c_{k}\right| +\left| c_{l}\right| \leq c_{0}\left(1+\cos \frac{\pi }{[
\frac{n-1}{k+l}]+2}\right)^{1/2}
\left(1+\cos \frac{\pi }{[\frac{n-1}{\left| k-l\right|
}]+2}\right)^{1/2},
\end{equation*}
for any distinct numbers $k$ and $l$ among $\{0, \ldots, n-1\}$.

\noi 

{\bf Organization of the paper.} We consider in the first Section
two classes of operator radii, called admissible and strongly
admissible radii. The operator norm and the numerical radius belong to
both classes as well as the more general radii $\om_{\rho}$ for $\rho
\leq 2$. We prove in Section 2 some
constrained and unconstrained von Neumann
inequalities for (strongly) admissible radii using the construction of
analytic
models of \cite{aem}. In Section 3 we prove some constrained von
Neumann inequalities for radii which are associated to some bundles of
operators ; these radii are not necessarily admissible. The
constraints in Section 2 are of algebraic type ($p(T) = 0$ or
$P(T^*,T) = 0$) while in Section 3 they are of the type $u(T) = 0$ for
a inner function $u$. Several
applications of the above general constrained von Neumann
inequalities are given in Section 4. Applications to bounds of
positive rational functions are presented in Section 5. In the last
section we discuss constrained von Neumann inequalities with
different type of constraints.

\section{Admissible and strongly admissible operator radii}

\noi 

{\bf Admissible operator radii.}
In this paragraph $w$ denotes a family of so-called operator
radii $w = \{w_{H}\}$, one for each separable Hilbert space under
consideration. An operator radius $w_{H}$ is a map from
$\LL(H)$ to $[0,+\infty]$. For $T \in \LL(H)$ we simply write $w(T)$
instead of the more correct $w_{H}(T)$ and say that $w$ is an
\emph{operator radius}, or simply a radius.

\begin{definition}
 A radius $w$ defined for all Hilbert space operators with values
 in $[0,+\infty]$ is called an
\emph{admissible radius} if it satisfies
\begin{itemize}
\item[(i)] ({\bf unitary invariance}) $w(U^*TU) = w(T)$ for each
unitary $U : K \to H$ and each $T \in \LL(H)$ ;
\item[(ii)] ({\bf isotonicity for restrictions}) if $T\in \LL(H)$ and
$E \subset H$ is invariant for $T$, then $w(T\mid E) \leq
w(T)$ ;
\item[(iii)] ({\bf ampliation}) if $T^{(\infty)}$ denotes the
countable orthogonal sum $T\oplus T \oplus \cdots$, then
$w(T^{(\infty)}) = w(T)$.
\end{itemize}
\end{definition}
The order on the extended interval
$[0,+\infty]$ uses the
usual conventions. In most examples we are looking for radii with
finite values.

\begin{remark} Suppose condition (i) holds. Then the ampliation
axiom (iii) is equivalent to
\begin{itemize}
\item[(iii')]  $w(T^{(n)}) = w(T)$, for $n$ finite or $n = \infty$,
\end{itemize}
where $T^{(n)}$ denotes the orthogonal sum of $n$ copies of $T$.
Indeed, there is a unitary equivalence between
$(T^{(n)})^{(\infty)}$ and $T^{(\infty)}$.
By \cite[Lemma 15]{fong/holbrook}, the ampliation condition is also
equivalent to $w( T \otimes I_{E}) = w(T)$, where $T\in \LL(H)$ and
$T \otimes I_{E}\in
\LL(H\otimes E)$. Here $I_{E}$ is the identity on $E$.
We refer to \cite{fong/holbrook} for other possible axioms
of operator norms and several examples.

Note also that half of condition (iii), namely
$w(T) \leq w(T^{(\infty)})$, is implied by conditions (i) and (ii).
\end{remark}

Let $\mathcal{F} = \{\mathcal{F}_{H}\}$
be a collection of Hilbert space operators, that is
for each considered
separable Hilbert space $H$,
$\mathcal{F}_{H} = \mathcal{F} \cap \LL(H)$ is a given
set.
\begin{definition}
Let $\mathcal{F}$ be a collection of Hilbert space operators. We say
that $\mathcal{F}$ is \emph{admissible} if it satisfies
\begin{itemize}
\item[(i)] ({\bf unitary invariance}) if $T \in \mathcal{F}\cap \LL(H)$
and $U : K \to H$ is unitary,
then $U^{\ast}TU \in \mathcal{F}\cap \LL(K)$ ;
\item[(ii)] ({\bf stability for restrictions}) if $T\in \mathcal{F}\cap
\LL(H)$
and $E \subset H$ is invariant for $T$, then $T\mid E \in \mathcal{F}$ ;
\item[(iii)] ({\bf ampliation}) if $T\in \mathcal{F}\cap \LL(H)$, then
$T^{(\infty)} \in \mathcal{F}$.
\end{itemize}
\end{definition}

\noi 

{\bf Radius associated to a collection of operators.} Let
$\mathcal{F}$ be a collection of Hilbert space operators.
Define the radius $w_{\mathcal{F}}$ associated to $\mathcal{F}$
by setting, for $T\in \LL(H)$,
$$w_{\mathcal{F}}(T) := \inf \{r > 0 : \frac{1}{r}T \in 
\mathcal{F}\cap \LL(H)\} .$$

\begin{prop}
    The radius associated to an admissible collection is an admissible
    radius.
\end{prop}

\begin{proof} Let
$\mathcal{F}$ be an admissible collection. In order to show the
unitary invariance of $w_{\mathcal{F}}$ let $T \in \LL(H)$  and
let $U : K \to H$ be a unitary operator. Fix $\ep > 0$. There exists $r =
r(\ep)$ such that $0 < r < w_{\mathcal{F}}(T) + \ep$ and
$\frac{1}{r}T \in \mathcal{F}\cap \LL(H)$. By the unitary invariance
of $\mathcal{F}$, we have $\frac{1}{r}U^{\ast}TU \in \mathcal{F}\cap
\LL(K)$. This shows that $w_{\mathcal{F}}(U^{\ast}TU) \leq r <
w_{\mathcal{F}}(T) + \ep$. Since $\ep > 0$ is arbitrary, we get
\begin{equation}
    w_{\mathcal{F}}(U^{\ast}TU) \leq w_{\mathcal{F}}(T) .
    \label{eq:leq}
\end{equation}
Replacing in
this inequality $T$ by $UTU^{\ast}$ we obtain $w_{\mathcal{F}}(T)
\leq w_{\mathcal{F}}(UTU^{\ast})$ ; replacing now $U$ by $U^{\ast}$ we
get
\begin{equation}
    w_{\mathcal{F}}(T) \leq w_{\mathcal{F}}(U^{\ast}TU) .
    \label{eq:geq}
\end{equation}
Using (\ref{eq:leq}) and (\ref{eq:geq}) we get that 
$w_{\mathcal{F}}$ is unitarily invariant.

The inequalities $w_{\mathcal{F}}(T\mid E) \leq
w_{\mathcal{F}}(T)$ and $w_{\mathcal{F}}(T^{(\infty)}) \leq
w_{\mathcal{F}}(T)$ can be proved as (\ref{eq:leq}) by using
the stability for restrictions and 
the ampliation axiom for $\mathcal{F}$, respectively.
Since $T$ is unitarily equivalent to a
restriction of $T^{(\infty)}$ to an invariant subspace, we also
obtain $w_{\mathcal{F}}(T) \leq w_{\mathcal{F}}(T^{(\infty)})$. Thus
$w_{\mathcal{F}}$ is admissible.
\end{proof}

In order to present some examples of admissible collections, we
introduce the following notation.
If $z$ is the variable in the complex plane $\C$, we denote by $P(\C)$ the
algebra of all complex polynomial functions in $\overline{z}$ and $z$.
If $T \in \LL(H)$  and $P\in P(\C)$, $P(\overline{z},z) =
\sum_{\alpha,\beta}c_{\alpha,\beta}\overline{z}^{\alpha}z^{\beta}$, we
set
$$P(T^*,T) = \sum_{\alpha,\beta}c_{\alpha,\beta} T^{\ast\alpha}T^{\beta} .$$
This is part of the so-called "hereditary functional calculus"
\cite{agler:ieot} which is briefly described in the next section.
We denote by $\sigma(T)$ the spectrum of an operator $T \in \LL(H)$.

\begin{theorem}\label{thm:adm}
    Let $\{P_{\lambda}\}_{\lambda \in \Lambda}$
    be a family of elements in
    $P(\C)$.

    (a) \quad Let $\mathcal{F}$ be the collection of operators defined by
    the following positivity conditions
    $$T \in \mathcal{F}\cap \LL(H) \mbox{ if and only if }
    P_{\lambda}(T^{\ast},T) \geq 0 \; (\lambda \in \Lambda).$$
    Then $\mathcal{F}$ is
    admissible.

    (b) \quad Define the
    collection $\mathcal{G}$ by
    $$T \in \mathcal{G}\cap \LL(H) \mbox{ if and only if }
     \sigma(T)  \subseteq \overline{\D}
     \mbox{ and }
     P_{\lambda}(T^{\ast},T) \geq 0 \; (\lambda \in \Lambda) .$$
    Then $\mathcal{G}$ is
    admissible.
\end{theorem}

\begin{proof} Let $U : K \to H$ be a unitary operator and let
$T\in \LL(H)$.
For each $\lambda\in \Lambda$ we have
$P_{\lambda}(U^{\ast }TU)=U^{\ast }P_{\lambda}(T)U$. Therefore
$\mathcal{F}$ is unitarily invariant.
If $E$ is an invariant subspace for $T$, then $(T\mid
E)^{\beta} = (T^{\beta})\mid E$ and $(T\mid E)^{\ast\alpha} =
P_{E}T^{\ast\alpha}\mid E$, where $P_{E}$ is the orthogonal projection
onto $E$. This shows that for each $\lambda\in \Lambda$ we have
$P_{\lambda}(T\mid E)=P_{E}P_{\lambda}(T)\mid E$, yielding the
stability to restrictions property.
The ampliation condition follows
from the equality $P_{\lambda}(T^{(\infty)}) = P_{\lambda}(T)^{(\infty)}$.

For the second part, note that the spectrum
satisfies $\sigma(U^{\ast }TU) = \sigma(T)$ and $\sigma(T^{(\infty)}) =
\sigma(T)$. Let now $T \in \LL(H)$ with $\sigma(T) \subseteq \overline{\D}$.
Let
$R = T\mid E \in \LL(E)$, where $E$ is an invariant subspace of $T$.
Thus the matrix of $T$ with respect to the decomposition $H = E
\oplus E^{\perp}$ has the form
$$T = \left(\begin{array}{cc}
            R & \ast  \\
            0 & \ast
        \end{array}\right) $$
and thus
$$T^n = \left(\begin{array}{cc}
            R^n & \ast  \\
            0 & \ast
        \end{array}\right) .$$
This implies $\|R^n\| \leq \|T^n\|$, so the spectral radius of $R$ is
at most one. This completes the proof.
\end{proof}
\begin{remark}
    Part (a) of the above
    Theorem also holds for more general positivity conditions,
    obtained by considering polynomials in $\overline{z}$ and $z$ with
matrix
    coefficients. We omit the details.
    Bounded collections satisfying such more general
    positivity conditions were characterized by J. Agler
    \cite{agler:models} as bounded collections which are closed with
    respect to direct sums, with respect to unital $C^{\ast}$-algebraic
    representations and stable for restrictions.
    We refer to \cite{agler:models} for the exact definition
    and for several examples of such collections.
\end{remark}

\noi 

{\bf Operators of class $C_{\rho}$.} The main examples of
operator radii we will use are the operator
radii
associated to the collection of
operators of class $C_{\rho}$.

Operators in the class
$C_{\rho}$ are
defined as operators having $\rho$-dilations :
$T \in \LL(H)$ is in $C_{\rho}$, $\rho > 0$, if there
exist a larger Hilbert space $K \supset H$ and a unitary operator
$U\in \LL(K)$ such that
$$T^nh = \rho P_{H}U^{n}h , \quad h \in H \; .$$
Contractions are operators of class $C_{1}$ and operators in $C_{2}$
coincides with
numerical radius contractions, that is operators
$T$ such that $\om_{2}(T) \leq 1$. We refer to \cite{SNF} and
\cite{racz} for more
information.

The operator radius $\om_{\rho}$ associated to the class
$C_{\rho}$ is then defined by
$$\om_{\rho}(T) = \inf \{ r : r > 0, \frac{1}{r}T \in C_{\rho}\}\; .$$
It is determined by the conditions that it is homogeneous
($\om_{\rho}(zT) = |z|\om_{\rho}(T)$ for all complex $z$) and that
$\om_{\rho}(T) \leq 1$ if and only if $T \in C_{\rho}$. Then $\om_{1}(T) =
\|T\|$ and $\om_{2}$ is the numerical radius. It can also be proved
that the limit of $\om_{\rho}(T)$ as $\rho\to \infty$ is the spectral
radius of $T$.

The radius $\om_{\rho}$ is a (Banach space)
norm if and only if $\rho \leq 2$. It is not an algebra norm ;
however, we always have  \cite{SNF}
$\om_{\rho}(T^n) \leq \om_{\rho}(T)^n$.

\begin{cor}
    The radius $\om_{\rho}$
    is admissible
    for any $\rho > 0$.
\end{cor}

\begin{proof} An operator $T$ is in the 
class $C_{\rho }$ \cite{SNF}, \cite{racz}
if and only if
$$
\left\| x\right\| ^{2}-2\left(1-\frac{1}{\rho}\right)
\mbox{ Re } \left[z<Tx\mid x>\right] +
\left(1-\frac{2}{\rho }\right) \left| z\right| ^{2}\left\|
Tx\right\| ^{2}\geq 0
$$
for every $x\in H$ and every $z\in \overline{\D}$.
Therefore it suffices to set
$$
P_{\lambda}(\overline{z},z) =
1-\left(1-\frac{1}{\rho}\right)\overline{\lambda}\overline{z} -
\left(1-\frac{1}{\rho}\right)\lambda z + \left(1-\frac{2}{\rho}
\right)|\lambda |^2\overline{z}z \quad  \quad (\lambda \in \overline{\D})
$$
in Theorem \ref{thm:adm}, Part (a).
\end{proof}

\noi 

{\bf Strongly admissible operator radii.}
The following definition gives a smaller class of admissible radii.

\begin{definition}
A radius $\nu$ defined for all Hilbert space operators with values
 in $[0,+\infty]$ is called a
\emph{strongly admissible radius} if it satisfies
\begin{itemize}
\item[(ii')] ({\bf isometry growth condition}) For any isometry $
V : K \to H$ and any $T\in B(H)$, we have $\nu
(V^{\ast }TV)\leq \nu (T)$.
\item[(iii)] ({\bf ampliation})
$\nu(T^{(\infty)}) = \nu(T)$ for every $T$.
\end{itemize}
\end{definition}
\begin{prop}
    An operator radius $\nu$ is strongly admissible if and only if it
    satisfies
    \begin{itemize}
    \item[(i)] ({\bf unitary invariance}) $\nu(U^*TU) = \nu(T)$ for each
unitary $U : K \to H$ and each $T \in \LL(H)$ ;
    \item[(ii')] ({\bf isotonicity for compressions}) If $T\in \LL
    (H)$, if $E$ is a
    closed subspace of $H$ and $R =
    P_{E}T\mid E$, then $\nu(R) \leq \nu (T)$.
    \item[(iii)] ({\bf ampliation}) We have
$\nu(T^{(\infty)}) = \nu(T)$.
\end{itemize}
In particular, each strongly admissible radius is admissible.
\end{prop}

\begin{proof} Suppose that $\nu$ is strongly admissible. Let $T \in
\LL(H)$ and let $U : K \to H$ be a unitary operator. Using the isometry growth
condition for the isometry $U$ we obtain $\nu(T_{1}) \leq \nu(T)$,
where $T_{1} =U^*TU$. The isometry growth
condition for $U^*$ yields $\nu(UT_{1}U^*) \leq \nu(T_{1})$. Therefore
$$
\nu(T) = \nu(UT_{1}U^*) \leq \nu(T_{1}) \leq \nu(T)
$$
showing the unitary invariance. The isotonicity for compressions is
obtained from $R = J^*TJ$, where $J : E \to H$ is the inclusion $Je =
e$.

For the converse implication, note that every isometry $V : K \to H$
can be written as $V = JU$, where $U : K \to V(K)$, $Uk = Vk$, is
unitary and $J : V(K) \to H$ is the inclusion map. Then
$$\nu(V^*TV) = \nu(U^*J^{*}TJU) = \nu(J^{*}TJ) = \nu(P_{V(K)}T\mid
V(K)) \leq  \nu(T).$$
The proof is now complete.
\end{proof}

A counterpart notion of \emph{strongly admissible collection} of operators
can be introduced as a collection which satisfies the unitary
invariance, the stability for
compressions and the ampliation properties. The radius associated to
a strongly admissible collection is strongly
admissible. We omit the details.

\begin{prop}
    Let $\rho > 0$. The radius $w_{\rho}$ is stongly admissible if and
    only if $\rho \leq 2$, if and only if $w_{\rho}$ is a norm.
\end{prop}

\begin{proof} Suppose $\rho \leq 2$.
Recall that $T \in C_{\rho}$ if and only if
$P_{\lambda}(T^{*},T) \geq 0$ for all $\lambda \in \overline{\D}$,
where
$$
P_{\lambda}(\overline{z},z) =
1-\left(1-\frac{1}{\rho}\right)\overline{\lambda}\overline{z} -
\left(1-\frac{1}{\rho}\right)\lambda z + \left(1-\frac{2}{\rho}
\right)|\lambda |^2\overline{z}z.
$$
The isometry growth condition for $w_{\rho }$ follows from the fact
that $V^{*}TV \in C_{\rho}$ whenever $\rho \leq 2$, $T\in C_{\rho}$
and $V : K\to H$ satisfies $V^*V = I_{K}$. Indeed, we have
\begin{eqnarray*}
    P_{\lambda}(V^{*}T^*V,V^{*}TV) & = & I-\left(1-\frac{1}{\rho
    }\right)\overline{\lambda}V^{*}T^{*}V -\left(1-\frac{1}{\rho
    }\right)\lambda V^{*}V  \\
     &  & + \left(1-\frac{2}{\rho }\right)|\lambda |^2V^{*}T^{*}VV^{*}TV \\
     & = & V^{*}P_{\lambda}(T^*,T)V +
     \left(\frac{2}{\rho}-1\right)|\lambda |^2V^{*}\left[ T^{*}
     (I-VV^*)T\right] V .
\end{eqnarray*}

Suppose now $\rho >2$ and let $b > 0$ be a fixed, arbitrary positive number.
Consider the following $2\times 2$ matrix
$$
T=\left(
\begin{array}{cc}
1 & b \\
0 & -1
\end{array}
\right) .
$$
We have \cite[Theorem 6]{ando/nishio}
$$
w_{\rho }(T)=\frac{1}{\rho }[\sqrt{\frac{\left| b\right| ^{2}}{4}+1}+\sqrt{
\frac{\left| b\right| ^{2}}{4}+1+\rho (\rho -2)}] \; ;
$$
in particular
$w_{2}(T)=[\left| b\right| ^{2}/4+1]^{1/2}.$
We can find a vector $e$ $\in \C^{2}$ such that $\left\| e\right\| =1$ and $
w_{2}(T)=\left| <Te\mid e>\right| $. Denote by $V$ the isometry from $\C$
into $\C^{2}$ defined by $V(z)=ze$. We have $V^{\ast }TV=<Te\mid e>e\otimes e$.
Therefore
$$
w_{\rho }(V^{\ast }TV)=\left| <Te\mid e>\right| =\sqrt{\left| b\right|
^{2}/4+1}.
$$
We have $w_{\rho }(V^{\ast }TV)>w_{\rho }(T)$ for any $\rho >2$.
It follows that $w_{\rho}$ is not a strongly admissible radius if $\rho
>2$. Recall \cite{SNF}
also that $w_{\rho }$ is a norm if and only if $\rho \leq
2$. 
\end{proof}
\begin{remark}
There are other interesting
examples of admisible and strongly admissible radii. For instance, if
$$W(T) = \{\langle Tx|x\rangle | : x \in H, \| x \| = 1 \}$$
denotes the numerical range of $T$, then the diameter of $W(T)$
$$\mbox{ diam }W(T) = \sup \{|\lambda - \mu | : \lambda , \mu \in
W(T)\}$$
is a
strongly admissible radius. Indeed (see for instance
\cite{gust/rao:book} for properties of the numerical range),
we have $W(U^*TU) = W(T)$, $W(P_{E}T\mid E) \subseteq W(T)$ and
$W(T^{(\infty)}) = W(T)$. Note also that the sum, or even convex
combinations, of (strongly) admissible radii are (strongly)
admissible. For instance, $T \to \|T\| + \mbox{ diam }W(T)$ is
strongly admissible.
\end{remark}

\section{(Constrained) von Neumann inequalities using
analytic models}
The existence of a model for contractions is a key result in Sz.-Nagy
and Foias dilation theory. In particular, a Hilbert space contraction
with spectrum contained in the open unit disc is unitarily equivalent
to a restriction of the backward shift of infinite multiplicity to an
invariant subspace. This implies easily inequality (\ref{vNineq2}) for
strict contractions. If $T$ is an arbitrary contraction, then, for
any real $r < 1$, inequality (\ref{vNineq2}) holds for the strict
contraction $rT$. Making $r\to 1$ we obtain (\ref{vNineq2}) for all
contractions.

We show in this section how the existence of a model implies at once von
Neumann and constrained von Neumann inequalities for different
admissible radii. We use the recent construction of analytic models
for $n$-tuples of operators due to Ambrozie, Engli\v{s} and M\"uller
\cite{aem}.

\noi 

{\bf Hilbert spaces associated to a domain.}
We recall the context of \cite{aem}, with some change of notation. We
refer to \cite{aem} and the references cited therein for more
information.

Let $D$ be a nonempty open domain in $\C^n$. Set $D^{\ast} = \{
\overline{z} : z \in D\}$. Let $\H$ be a
$D$-\emph{space}, that is $\H$ is a Hilbert space of functions
analytic on $D$ such that
\begin{itemize}
    \item[(a)] $\H$ is invariant under the operators $Z_{j}$, $j = 1,
    \ldots n$, of multiplication by the coordinate functions,
    $$(Z_{j}f)(z) = z_{j}f(z) \quad ; \quad f\in \H \; ,\;  z = (z_{1},
    \ldots , z_{n})\in D.$$

    \item[(b)] For each $z\in D$, the evaluation functional $f\to f(z)$
    is continuous on $\H$.

    \item[(c)] $C(w,z) \neq 0$ for all $z \in D$ and $w \in D^{\ast}$.
\end{itemize}
Here $C(w,z)$ is the reproducing kernel of $\H$, that is
$C(w,z) = C_{\overline{w}}(z)$, for $z \in D$ and $w \in D^{\ast}$,
where $C_{\zeta}$ is a function in $\H$ such that $f(\zeta) = \langle
f\mid C_{\zeta}\rangle$, $f\in \H$
(we use (b) and the Riesz representation theorem).

Let $H$ be a Hilbert space. Denote by $\H\otimes H$ the completed
Hilbertian tensor product. Consider the multiplication operators
$M_{z_{j}}$ on $\H\otimes H$ defined by
$$M_{z_{j}} = Z_{j}\otimes I_{H} \quad ; \quad j = 1, \ldots , n .$$
Set
$$Z = (Z_{1},\ldots, Z_{n}) \in \LL (\H)^n \quad ; \quad
M_{z} = (M_{z_{1}}, \ldots  , M_{z_{n}}) \in \LL (\H\otimes H)^n.$$

Let $T = (T_{1}, \ldots, T_{n})$ be a commuting tuple of operators.
Denote by $\sigma (T)$ the Taylor spectrum of $T$, and let
$$M_{T} = (L_{T_{1}^{\ast}}, \ldots L_{T_{n}^{\ast}}, R_{T_{1}}, \ldots
R_{T_{n}}) .$$
Here $L_{A}(X) = AX$ and $R_{A}(X) = XA$ are the left and right
multiplication operators by $A$ on $\LL (H)$. Let $F$ be a analytic
function on a neighborhood of $\sigma (M_{T})$. Define
$F(T^{\ast},T)
\in \LL(H)$ by $F(T^{\ast},T) = F(M_{T})(I)$.

If $z = (z_{1}, \ldots, z_{n})$ is the variable
in complex Euclidean space $\C^n$, we denote by
$P(\C^n)$ the
algebra of all complex polynomial functions in
$\overline{z}_{1}, \ldots, \overline{z}_{n}$, $z_{1}, \ldots, z_{n}$.
If $F(w,z) =
w^{\alpha}z^{\beta}$, then $F(T^{*},T) = T^{*\alpha}T^{\beta} =
P(T^{*},T)$ for $P(\overline{z},z) = F(\overline{z},z)
\in P(\C^n)$. We use the usual notation $T^{\beta} =
T_{1}^{\beta_{1}}\cdots T_{n}^{\beta_{n}}$ for $\beta = (\beta_{1},
\ldots \beta_{n}) \in \N^n$ and the like.
Note that this differs slightly from \cite{aem} where
$T^{\ast}$ is written on the right.

\noi 

{\bf Axiom (AEM).}
We will sometimes suppose that $\H$ satisfies \emph{Axiom (AEM)}, that is $\H$
is a $D$-space
such that the polynomials are dense
in $\H$ and $\frac{1}{C}$ is a polynomial.
Let $(\psi_{k})$ be a fixed
orthonormal basis for $\H$ consisting of polynomials such that any
finite polynomial is a finite linear combination of $\psi_{k}$. Set
$$f_{m}(w,z) = \sum_{k=m}^{+\infty}
\overline{\psi_{k}(\overline{w})}\frac{1}{C}(w,z)\psi_{k}(z) .$$

When $D = \D$ is the open unit disk and $\H$ is the Hardy space $H^2 =
H^2(\D)$
of the unit disk, then $C(w,z) = (1-wz)^{-1}$ and $M_{z}^*$ is the
backward shift of multiplicity $\dim H$. In this case $\H$ satisfies
axiom (AEM) with $\psi_{k}(z) = z^k$
and $f_{m}(w,z) = w^mz^m$. We refer to \cite{aem} for
other examples.

\noi 

{\bf Unconstrained von Neumann inequalities for operator radii.}
We use notation as above.

\begin{theorem}\label{thm:unco}
    Let $T = (T_{1}, \ldots, T_{n}) \in \LL(H)$
    be an $n$-tuple of commuting operators. Suppose $T$ and $\H$
    satisfy one of the following two conditions
    \begin{itemize}
    \item[(i)] $\H$ is a $D$-space, $\sigma (T) \subset D$ and
    $\frac{1}{C}(T^{\ast},T) \geq 0$ ;
    \item[(ii)] $\H$ is a $D$-space satisfying Axiom (AEM),
    $\frac{1}{C}(T^{\ast},T) \geq 0$ and
    $$\lim_{m}f_{m}((T^{\ast},T))h = 0$$ 
    for every $h\in H$.
\end{itemize}
    Let $p(z) = \sum_{\beta\in \Z_{+}^n}c_{\beta}z^\beta$ be
    a fixed polynomial in the variable $z\in \C^n$ and let $P(w,z) =
    \sum_{\alpha,\beta\in \Z_{+}^n}c_{\alpha,\beta}
    w^\alpha z^\beta$ be a fixed polynomial
    in two variables. If $\om$
is an admissible radius, then
$$\om (p(T)) \leq \om (p(Z^*)) \; ;$$
if $\nu$ is strongly admissible, then
$$\nu (P(T^{\ast},T)) \leq \nu (P(Z,Z^*)) .$$
\end{theorem}

\begin{proof}
Suppose $T$ satisfies (i) or (ii).
In either case, using \cite[Corollary 7,Corollary 15]{aem},
there is an isometry $V : H \to \H\otimes H$ such that $VT_{j} =
M_{z_{j}}^*V$ for $j = 1, \ldots, n$. Note again that some care has to be
taken when using the results of \cite{aem} because of the change of
notation.
This implies
$$Vp(T) = p(M_{z}^{\ast})V.$$
In particular, $VH$ is invariant under
$M_{z}^{\ast}$ and $T^{\beta}$ is unitarily
equivalent to the restriction of $M_{z}^{\ast\beta}$ to the invariant
subspace $VH$. Since $\om$ is admissible, we have
$$\om (p(T)) = \om (p(M_{z}^{\ast})\mid VH) \leq \om
(p(M_{z}^{\ast}).$$
Using the ampliation axiom for $\om$ and the fact that
$M_{z}^* = Z^*\otimes I_{H}$, we obtain
$$\om (p(T)) \leq \om (p(Z^*)) .$$

For the second part of the theorem, note that with respect to
the decomposition $\H\otimes H = VH \oplus VH^{\perp}$, we can write
$$M_{z}^{\ast \beta} =
\left(\begin{array}{cc}
            T^{\beta} & \ast \\
            0 & \ast
        \end{array}\right) \quad \mbox{ and } \quad M_{z}^{\alpha} =
        \left(\begin{array}{cc}
            T^{*\alpha} & 0 \\
            \ast & \ast
        \end{array}\right) .$$
This shows that
$$T^{*\alpha}T^{\beta} = P_{VH}M_{z}^{\alpha}M_{z}^{\ast \beta}|VH .$$
Since $\nu$ is strongly admissible, we have
$$\nu (P(T^*,T)) \leq \nu (P(M_{z},M_{z}^{\ast})) = \nu (P(Z,Z^*)) .$$
The proof is complete. 
\end{proof}

\noi 

{\bf Constrained von Neumann inequalities.}
We start with a constrained von Neumann inequality for admissible radii.

\begin{theorem}\label{thm:24}
    Let $D$ be an open domain in $\C^n$. Suppose Hilbert space $\H$ of
    functions analytic on $D$ and an $n$-tuple of operators $T$
    satisfy one of the two conditions (i) and (ii)
    in Theorem \ref{thm:unco}.
    Let $p$ and $q$ be one variable polynomials in $n$ variables
    and suppose that $q(T) = 0$. If $\om$
is an admissible radius, then
$$\om (p(T)) \leq \om (p(Z^*|\Ker q(Z^*))) .$$
\end{theorem}

\begin{proof} We use the notation of (the proof of) Theorem
\ref{thm:unco}.
Recall that $T^{\beta}$ is unitarily
equivalent to the restriction of $M_{z}^{\ast\beta}$ to the invariant
subspace $VH$ and $Vp(T) = p(M_{z}^{\ast})V$. Since $q(T) = 0$, we have
$$0 = Vq(T)h = q(M_{z}^{\ast})Vh$$
for any $h \in H$. This shows that $VH \subseteq \Ker q(M_{z}^{\ast})$.
Therefore
$$\om (p(T)) = \om (p(M_{z}^{\ast})\mid VH) \leq \om
(p(M_{z}^{\ast})\mid \Ker q(M_{z}^{\ast})).$$
By the ampliation axiom and the equality 
$$M_{z}^*\mid \Ker q(M_{z}^{\ast})
 = [Z^*\mid \Ker q(Z^*)]\otimes I_{H},$$ 
we get
$$\om (p(M_{z}^{\ast})\mid \Ker q(M_{z}^{\ast})) = \om(p(Z^*\mid \Ker
q(Z^*))).$$ 
This completes the proof.
\end{proof}

In some applications it is possible to avoid the
hypothesis 
$$\lim_{m}f_{m}((T^{\ast},T))h = 0$$ 
in
condition (ii) in Theorem \ref{thm:unco}. We refer to Corollary
\ref{cor:41} and Corollary \ref{cor:46}
for examples of results of this type.

The following result is a constrained von Neumann inequality for
strongly admissible radii. Recall 
that $E\subseteq H$ is said to be \emph{invariant} 
for the $n$-tuple $(T_{1},\ldots,
T_{n}) \in \LL(H)$ if $T_{j}E \subseteq E$ for each $j$.

\begin{theorem}\label{thm:25}
    Let $D$ be an open domain in $\C^n$. Suppose Hilbert space $\H$ of
    functions analytic on $D$ and an $n$-tuple of operators $T$
    satisfy one of the two conditions (i) and (ii)
    in Theorem \ref{thm:unco}. Suppose also that each operator $Z_{j}\in
    \LL(\H)$
    is an isometry.
    Let $P$ and $Q$ be two elements of $P(\C^n)$
    and suppose that $Q(T^*,T) = 0$. There exists an invariant subspace
    $E$ for $Z^* \in \LL(\H)^n$ such that,
    for each strongly admissible radius $\nu$,
$$\nu(P(T^*,T)) \leq \nu(P(Z_{E},Z_{E}^*)),$$
where $Z_{E}$ is defined by setting $Z_{E}^* := Z^* \mid E$. If $n = 1$, 
if $Q$ is
of degree less or equal than $d$ and if
$Q(e^{-it},e^{it})\neq 0$ for some $t \in \R$, then the dimension of $E$ is
less or equal than $2d$.
\end{theorem}

\begin{proof} Let $Q$ be a polynomial in $P(\C^n)$ of degree at
most $d$, that is, the maximal power at which each $w_{j}$ and $z_{j}$
occurs is at most $d$. Recall from the proof of Theorem \ref{thm:unco}
that
$$T^{*\alpha}T^{\beta} = P_{VH}M_{z}^{\alpha}M_{z}^{\ast \beta}|VH$$
and thus $P(T^*,T) = J^*P(M_{z},M_{z}^*)J$, where $J$ denotes the
inclusion $J : VH \to \H\otimes H$. The same equality, using the fact
that $Q(T^*,T) = 0$, implies that the subspace $V(H)$ is contained
in $\Ker Q(M_{z},M_{z}^*)$. Since each $Z_{j}$ and thus each
$M_{z_{j}}$ is an isometry, we get that $\Ker Q(M_{z},M_{z}^*)$ is
included in $E_{0}=\Ker
(M_{z}^{\ast d}Q(M_{z},M_{z}^*))$
which is invariant by $M_{z}^*$. Denote
$$E = \Ker (Z^{\ast d}Q(Z,Z^*))$$
which is invariant under $Z^*$. Then,
using the properties of the strongly admissible radius $\nu $, we obtain
\begin{eqnarray*}
\nu (P(T^{\ast },T)) &=&\nu (J^{\ast }P(M_{z},M_{z}^*)J) =
\nu (J^{\ast }P_{E_{0}}P(M_{z},M_{z}^*)P_{E_{0}}J) \\
&\leq &\nu (P_{E_{0}}P(M_{z},M_{z}^*)P_{E_{0}}) = \nu
(P_{E_{0}}[P(Z,Z^*)\otimes I]P_{E_{0}}) \\
& = &\nu (P_{E}(p(Z,Z^{\ast })P_{E}\otimes I)=\nu (P_{E}(p(Z,Z^{\ast })P_{E})
\\
&=&\nu (p(Z_{E},Z_{E}^{\ast }))
\end{eqnarray*}
where $P_{E_{0}}$ and $P_{E}$ are the orthogonal projections onto
$E_{0}$, respectivelly $E$. 

Finally, if $n = 1$ and if $Q$ is of degree less or equal than $d$, then 
$Z^{\ast d}Q(Z,Z^*)$ is a polynomial in $Z^{\ast}$ of degree less or equal to 
$2d$. Thus $E$, the kernel of $Z^{\ast d}Q(Z,Z^*)$, 
is a subspace of dimension 
no greater than $2d$, unless $Z^{\ast d}Q(Z,Z^*)$ is the null operator. This 
occurs if and only if $e^{ids}Q(e^{-is},e^{is}) = 0$ for every $s\in \R$. The 
last equality is impossible if $Q(e^{-it},e^{it})\neq 0$ for some $t \in \R$.
The proof is complete.
\end{proof}

\section{Inequalities for radii associated to bundles of operators}
Constrained von Neumann inequalities for some
operator radii which are not necessarily admissible are obtained in
this section. The method also gives a different proof of constrained
von Neumann inequalities for the radii $\om_{\rho}$.

\noi 

{\bf Notation.}
We denote by $D(\alpha,r)$ the open disc of radius $r$ and center
$\alpha$.
Let $\T$ be the boundary of $\D = D(0,1)$. The spaces $L^p = L^p(\T)$, $1\leq
p \leq \infty$, are the usual Lebesgue function spaces relative to
normalized Lebesgue measure on $\T$. The spaces $H^p = H^p(\T)$, $1\leq
p \leq \infty$, are the usual Hardy spaces. Denote
$$H^1_{0}= \{f \in L^1 : \int_{0}^{2\pi}f(e^{it})e^{int}\, dt = 0,
n = 0,1,\ldots \}.$$

For a given inner function $u$, denote $H(u) = H^2 \ominus uH^2$
and consider the operator
$S(u) \in \LL(H(u))$ defined by
$$S(u) = P_{H(u)}Z\mid H(u)\; .$$
Recall that $Z$ is the operator of multiplication by $z = e^{i\theta}$
on $H^2$. A proof that $S(u)$ and the extremal operator
$S^{\ast }\mid \Ker (u(S)^{\ast })$ are unitarily
equivalent follows from the fact that they have the
same characteristic function \cite{SNF} ; a direct proof can be found
in \cite{ptak:casopis}.

If $T\in \LL(H)$ is an absolutely continuous contraction, then, for
any $x,y \in H$, there exists a function $x\overset{T}{.}y\in L^{1}$
with the $n$th Fourier coefficient given by $\langle T^{*n}x\mid
y\rangle$ if $n\geq 0$ and $\langle T^{-n}x\mid
y\rangle$ if $n < 0$.

Let $T$ be an operator whose spectrum is included in the closed unit disc.
Consider the operator kernel $K_{\alpha}(T)$ defined by
$$
K_{\alpha }(T)=(I-\overline{\alpha }T)^{-1}+(I-\alpha T^{\ast
})^{-1}-I;\quad |\alpha | < 1.
$$
For an absolutely continuous contraction $T$, 
$<K_{r\exp(it)}(T)x\mid y>$
converges almost everywhere to $x\overset{T
}{.}y$ when $r$\ goes to $1$.

Recall that a contraction $T \in \LL(H)$ is said \cite{SNF}
to be of class $C_{0}$ if $T$ is c.n.u. and there is a nonzero
function $f$ in $H^{\infty}$ such that $f(T) = 0$. Then there is a
unique (up to a constant factor of modulus one) nonconstant inner
function $u$, called the \emph{minimal function} of $T$, such that $u(T) =
0$. The minimal function of $S(u)$ is $u$.

\noi 

{\bf Bundles of selfadjoint operators and associated radii.} Recall
the following result. Let $\rho > 0$. An operator $T\in\LL(H)$
whose spectrum is included
in the closed unit disc is in $C_{\rho }$
if and only if \cite{CaFa} $K_{\alpha }(T)+\rho I\geq I$
for any $\alpha \in \D$.

\begin{definition}
    Suppose a collection $\mathcal{R}$ of \emph{bundles} of self-adjoint
    operators is given,
    that is, for each separable Hilbert space $H$ there is a map
    $$\mathcal{R}_{H} : \D\times \LL(H) \ni (\alpha,A) \to
    R_{\alpha}(A) \in \LL(H)$$
    with $R_{\alpha}(A) = R_{\alpha}(A)^*$. The collection
    $\mathcal{K} = \mathcal{K}_{\mathcal{R}}$ associated to
$\mathcal{R}$ is
    defined by setting
    $$A\in \mathcal{K}\cap \LL(H) \mbox{ if and only if }
    \sigma(A) \subseteq \overline{\D} \mbox{ and } 
    K_{\alpha}(A) + R_{\alpha}(A) \geq I \;  \; (\alpha \in \D).$$
\end{definition}
The operator radius associated to $\mathcal{K} = \mathcal{K}_{\mathcal{R}}$
is then
$$\om_{\mathcal{K}}(A) = \inf \{r > 0 : \frac{1}{r}A \in
\mathcal{K}_{\mathcal{R}}\}.$$

\begin{example} \firsta Let $\rho > 0$. For the bundle $\mathcal{R}$ given by
$R(\alpha )=\rho I$, the class $\mathcal{K}_{\mathcal{R}}$ coincides with
the class $C_{\rho }$. 

\condb Let $A$ be an positive invertible operator and set $R(\alpha ) = A$.
Then
the associated collection $\mathcal{K}_{\mathcal{R}}$ coincides with
the class $C_{A}$ introduced by H. Langer (cf. \cite[p. 54]{SNF}). 

\condc Let $\varphi $ be a function in $H^{\infty }(\D)$ and let $v
: \LL(H) \to [0,+\infty[$ a map which satisfies
$$
v(T)\left\| T\right\| \leq 1\text{ for all }T\in B(H).
$$
Consider the bundle $\mathcal{R}$ by setting
$$
R_{\alpha }(T)=\varphi (\alpha v(T)T) + \varphi (\alpha v(T)T)^{\ast}.
$$
The associated radius to the collection $\mathcal{K}_{\mathcal{R}}$
is not necessarily admissible.
\end{example}

\begin{theorem}\label{thm:62} Let $\mathcal{R}$ be a bundle such that
$$\D \ni \alpha\to R_{\alpha}(A) \in \LL(H)$$
is harmonic in $\D$ for each $A \in \LL(H)$.
Let $\mathcal{K} = \mathcal{K}_{\mathcal{R}}$ be the collection of
operators associated to the bundle
$\mathcal{R}$.
Let $T$ be a contraction of class
$C_{0}$ with $u(T) = 0$, $u$ an inner function, and let $f\in A(D)$.
Assume that for any $\alpha \in D$ there exist a function
$
g_{\alpha }
$
such that for any $r>0$
$$
R_{\alpha }(f(T)/r) = g_{\alpha }(f(T)/r)+g_{\alpha }(f(T)/r)^{\ast } 
$$
and
$$
R_{\alpha }(f(S(u))/r) = g_{\alpha }(f(S(u))/r)+g_{\alpha }(f(S(u))/r)^{\ast} .
$$
Then we have
$$
\om_{\mathcal{K}}(f(T))\leq \om_{\mathcal{K}}(f(S(u)).
$$
\end{theorem}
Recall that $S(u)$ is unitarily equivalent to the extremal
operator $S^{\ast }\mid \Ker (u(S)^{\ast })$.

For the proof of Theorem \ref{thm:62},
we need the following lemma which will be also used in Section
\ref{section:5}.

\begin{lem}\label{lem:63}
Let $u$ be a inner function
    and let $f$ be a positive function in the subspace $\overline{u}
    H_{0}^{1}$ of $L^{1}(\T )$. Then there exists a function $h$ in $
    H^{2}\ominus uH^{2}$ such that $f=\left| h\right| ^{2}$.
\end{lem}

\begin{proof} Since $f$ $\in \overline{u}H_{0}^{1}$ we have $f=
\overline{u}f_{1}$, with $f_{1}\in H_{0}^{1}$. Then $\log
\left| f\right| = \log \left| f_{1}\right| $ is Lebesgue integrable.
According to theorem of Hoffman \cite{Hof} there exists an outer
function $g$ in $H^{2}$ such that $f=\left| g\right| ^{2}$. Denote by
$E = H(u)$
the orthogonal in $H^{2}$ of the subspace $uH^{2}$ and write $
g=g_{1}+ug_{2}$ with respect to the orthogonal decomposition $H^{2}=E\oplus
uH^{2}$. We have $g_{1}\neq 0$ since $g$ is an outer function.
Using the fact that $g_{1}\in E$, we obtain
$$
<u\overline{g_{1}}\mid \overline{h}>=\int_{0}^{2\pi }u(e^{it})\overline{
g_{1}(e^{it})}h(e^{it})dm(t)=\overline{<g_{1}\mid uh>}=0,
$$
for all functions $h$ in $H^{2}$.
Using the theorem of F. and M. Riesz \cite{Hof} we get
\begin{equation}
    u\overline{g_{1}}\in H_{0}^{2}.
    \label{eq:341}
\end{equation}
On the other hand, we have
\begin{eqnarray*}
uf=u\left| g\right|^{2} & = & u\left| g_1 + ug_2\right|^{2}\\
                        & = & u \left( g_1 + ug_2\right) 
            \left( \overline{g_1 + ug_2}\right)\\
      & = & u\left| g_{1}\right| ^{2}+u\left| g_{2}\right|
^{2}+g_{1}\overline{g_{2}}+u^{2}\overline{g_{1}}g_{2}.
\end{eqnarray*}
Therefore
$$
g\overline{g_{2}} = (g_1 + ug_2)\overline{g_{2}} 
= u\left| g_{2}\right| ^{2}+g_{1}\overline{g_{2}}=uf-u\left|
g_{1}\right| ^{2}-u^{2}\overline{g_{1}}g_{2}.
$$
Since $f\in \overline{u}H_{0}^{1}$ and using (\ref{eq:341}),
we see that the three last
terms belong to $H_{0}^{1}$. Hence $g\overline{g_{2}}\in
H_{0}^{1}$ and for any polynomial $p$ we have
$$
<pg\mid g_{2}>=\int_{0}^{2\pi }p(e^{it})g(e^{it})\overline{g_{2}(e^{it})}
dm(t)=0.
$$
Since $g$ is an outer function, it follows that $g_{2}=0$. The proof of
the lemma is now complete. 
\end{proof}

\begin{proof}[Proof of Theorem \ref{thm:62}]
By the canonical factorization theorem, $u$ \ can be decomposed
as
$$
u(z)=B(z)\exp [-\int_{0}^{2\pi }\frac{e^{i\theta }+z}{e^{i\theta }-z}d\mu
(\theta )]
$$
where $B$ is a Blaschke product and $\mu $ is a positive measure on $\partial D
$ which is singular with respect 
to the Lebesgue measure. Using the spectral
mapping theorem of a $C_{0}$ operator, we have
$$
\sigma (T)\subseteq \overline{B^{-1}\{0\}}\cup Supp(\mu )=\sigma (S(u)),
$$
where $Supp(\mu )$ is the support of $\mu $.

Let $f$ be a non-identically zero function in $A(\D)$. Using the spectral 
mapping 
theorem, we get
\begin{equation}
    \sigma (f(T))=f(\sigma (T))\subseteq f(\sigma (S(u)))=\sigma
    (f(S(u))).
    \label{eq:331}
\end{equation}
Fix $r > \om_{\mathcal{K}}(f(S(u)))$ and let 
$\alpha \in D(0,1/(r\left\| f\right\|_{\infty }))$, $\alpha \neq 0$. We
deduce from (\ref{eq:331}) that $r/\alpha$  
belongs to the resolvent of $T$.
Therefore, for any $x\in H$ and every $\alpha$ in
$D(0,1/(r\left\| f\right\|
_{\infty }))$, we can write
\begin{eqnarray*}
&<&[K_{\alpha }(\frac{f(T)}{r})+R_{\alpha }(\frac{f(T)}{r})-I]x\mid x> \\
&=&<\left[(I-\alpha \frac{f(T)^{\ast }}{r})^{-1}+(I-\overline{\alpha 
}\frac{f(T)}{
r})^{-1}-2I+g_{\alpha }(\frac{f(T)}{r})\right.\\ 
& & + \left. g_{\alpha }(\frac{f(T)}{r})^{\ast
}\right]x\mid x>.
\end{eqnarray*}
Recall that for any absolutely continuous
contraction $T$ and for any $x,y\in H$, the function $<K_{r,t}(T)x\mid y>$
converge almost everywhere
to a function $x\overset{T
}{.}y\in L^{1}(\partial D)$ when $r$\ goes to $1$.
Since $T$ is a $C_{0}$ contraction such that $u(T)=0$,
it follows \cite[Lemma 5.2]{CCC}
that $x\overset{T}{.}x\in \overline{u}H_{0}^{1}$. From
Lemma \ref{lem:63} we get the
existence of a function $h$ in $H^{2}\ominus uH^{2}=E$ such that $x\overset{T
}{.}x(e^{it})=\left| h(e^{it})\right| ^{2}$. We obtain
\begin{eqnarray*}
&<&[K_{\alpha }(\frac{f(T)}{r})+R_{\alpha }(\frac{f(T)}{r})-I]x\mid x> \\
&=&\int_{0}^{2\pi }\left[\frac{1}{1-\overline{\alpha
}\frac{f(e^{it})}{r})}+\frac{1}{
1-\alpha \overline{\frac{f(e^{it})}{r}}}-2 \right.\\ 
& & \left. + g_{\alpha
}(\frac{f(e^{it})}{r}) + \overline{g_{\alpha
}(\frac{f(e^{it})}{r})}\right]x\overset{T}{.}x(e^{it})dm(t) \\
&=&\int_{0}^{2\pi }\left[\frac{1}{1-\overline{\alpha
}\frac{f(e^{it})}{r})}+\frac{1}{
1-\alpha \overline{\frac{f(e^{it})}{r}}}-2 \right.\\  
& & \left. + g_{\alpha
}(\frac{f(e^{it})}{r}) + \overline{g_{\alpha
}(\frac{f(e^{it})}{r})}\right]\left| h(e^{it})\right| ^{2}dm(t) \\
&=&<\left[(I-\alpha \frac{f(S_{u})^{\ast }}{r})^{-1} + 
(I-\overline{\alpha }\frac{
f(S_{u})}{r})^{-1}-2I \right.\\ 
& & \left. + g_{\alpha }(\frac{f(S_{u})}{r})+ 
 g_{\alpha
}(\frac{f(S_{u})}{r})^{\ast }\right]h\mid h>
\\
&=&<\left[ K_{\alpha }(\frac{f(S_{u})}{r})+R_{\alpha
}(\frac{f(S_{u})}{r})-I\right] h\mid h>.
\end{eqnarray*}
Since both sides of the previous equalities are harmonic inside the unit disc
(with respect to the variable $\alpha $) and coincide inside
the disc $D(0,1/(r\left\| f\right\|_{\infty }))$,
we get that for any $\alpha $ in the unit disc
\begin{equation}\label{eq:33}
<\left[K_{\alpha }(\frac{f(T)}{r}) + 
R_{\alpha }(\frac{f(T)}{r})-I\right]x\mid x>
\end{equation}
and
\begin{equation}
<\left[ K_{\alpha }(\frac{f(S_{u})}{r}
)+R_{\alpha }(\frac{f(S_{u})}{r})-I\right] h\mid h>
\end{equation}
coincide. 
As $r>\om_{\mathcal{K}}(f(S_{u}))$, we get the positivity of
(\ref{eq:33}). We obtain $r>\om_{\mathcal{K}}(f(T))$ and the proof
is now complete. 
\end{proof}
\section{Applications of the previous results.}
We show in this section how the above constrained von Neumann
inequalities can be applied in a variety of situations. We are not
always looking for the most possible general inequalities.

\noi 

{\bf Applications of Theorem \ref{thm:24}.} We denote by $[ x ]$
the integer part of $x$, that is the least integer no greater than $x$.

\begin{cor}\label{cor:41}
Let $n\geq 2$. Let $T\in \LL(H)$ be a contraction such that $T^n =
0$. Then, for each $\rho > 0$ and each analytic polynomial $p$, we
have
$$\om_{\rho}(p(T)) \leq \om_{\rho}(p(S_{n}^*)).$$
In particular, for any $m$ we have
$$\om_{2}(T^m) \leq \cos \frac{\pi}{k(m,n)+2}, \quad k(m,n) :=
[\frac{n-1}{m}].$$
\end{cor}

\begin{proof} Let $r < 1$ be a positive number. Let $\H$ be $H^2$
in Theorem  \ref{thm:24}. Then
Theorem \ref{thm:24}, (ii),
applied with $rT$ instead of $T$, $q(z) = z^n$ and $\om = \om_{\rho}$,
gives
$$\om_{\rho}(p(rT)) \leq \om_{\rho}(p(S_{n}^*)).$$
Make now $r$ tends to $1$.

For the proof of the last part note that
a majorant of the left-hand side will
be $\om_{2}(S_{n}^{*m})$. But $S_{n}^{*m}$ is unitarily equivalent to an
orthogonal sum of shifts of smaller dimension, the largest dimension
being $k(m,n)+1$. Therefore
$\om_{2}(S_{n}^{*m}) = \om_{2}(S_{k(m,n)+1}^*)$ is equal
to $\cos \frac{\pi}{k(m,n)+2}$. The same computation follows
from \cite[page 120]{gust/rao:book}.
\end{proof}

\begin{remark}
The inequality $\om_{2}(T^m) \leq \cos (\pi/(k(m,n)+2))$
can be deduced from the inequality (\ref{HHineq}) of
Haagerup and de la Harpe. Indeed, $k = k(m,n) = [\frac{n-1}{m}]$ implies that
$mk+m > n-1$ and thus $(T^m)^{k+1} = 0$. Apply the Haagerup-de la
Harpe inequality for $T^m$.\label{rem:42}
\end{remark}

In the general case, if $p(z) = a_{0} + a_{1}z + \cdots +
a_{n-1}z^{n-1}$
is a polynomial of degree less or equal than  $n-1$, then $p(S_{n}^*)$
is the following triangular Toeplitz matrix
$$p(S_{n}^*) = \left(\begin{array}{ccccc}
            a_{0} & a_{1} & a_{2} & \cd & a_{n-1} \\
             & a_{0} & a_{1} & \cd & a_{n-2} \\
              & & a_{0} & \cd & a_{n-3} \\
             & & & \ddots & \vdots \\
            & & &  & a_{0}
        \end{array}\right). $$
Recall that we have the following reciprocity law of Ando and Nishio:
$$\omega _{\rho }(T)=(\frac{2}{\rho }-1)\omega _{2-\rho }(T) .$$
This shows that computations of $\omega _{\rho }$ for $0<\rho <1$
follows from computations for $1<\rho <2$.
Using interpolation properties of $\omega _{\rho }$
(see [FH, p.296]),
the law of Ando and Nishio, Corollary \ref{cor:41} and a result from
\cite[Lemma
5]{eiermann}
concerning the numerical range of Toeplitz matrices we get the next result.
\begin{cor}
Let $n\geq 2$. Let $T\in \LL(H)$ be a contraction such that $T^n =
0$. Let $p(z) = a_{0} + a_{1}z + \cdots +
a_{n-1}z^{n-1}$ be an analytic polynomial of degree at most $n-1$. 
Let $\theta \in \R$. We
have
$$
\omega _{\rho }(p(T))\leq \left(\frac{2}{\rho }-1\right)\left\| p\right\| 
_{\infty
}^{\rho }\left[ \inf_{\theta \in \R}\sup \{\left| p(\zeta )\right| :\zeta \in
\C,\zeta ^{2n-1}=e^{i\theta }\}\right] ^{1-\rho }$$ 
if $\rho \in ]0,1]$, and
$$
\omega _{\rho }(p(T))\leq \left\| p\right\| _{\infty }^{2-\rho }\left[
\inf_{\theta \in \R}\sup \{\left| p(\zeta )\right| :\zeta \in \C,\zeta
^{2n-1}=e^{i\theta }\}\right] ^{\rho -1}$$ 
if $\rho \in [1,2]$.
\end{cor}
We refer to the proof of Theorem \ref{thm:harm} for a better estimate
of $\om_{2}(p(S_{n}^*))$ for polynomials of the form
$p(z) = z^k + e^{i\gamma}z^l$ ; this yields (Theorem \ref{thm:662}) an
estimate for $\om_{2}(T^k+T^l)$.

If the polynomial $q$ of degree $d$ is given by
$$q(z) =    b_{0} + b_{1}z + \cdots +
b_{d}z^{d},$$
then $\Ker(q(S^*))$ consists of all sequences $(y_{r}) \in \ell_{2}$
satisfying
$$b_{d}y_{r+d} + b_{d-1}y_{r+d-1} + \cdots b_{0}y_{r} = 0$$
for $r = 0, 1, 2, \ldots$. This linear recurrence has a
$d$-dimensional solution space and if all the zeros of $q$ have
modulus less than one then all solutions lie in $\ell_{2}$. In this
case $\Ker(q(S^*))$ has dimension $d$.

We refer to \cite{ptak/young} for the matrix of  $S^*\mid\Ker(q(S^*))$
with respect to some orthonormal basis of  $\Ker(q(S^*))$ and,
for instance, to \cite{gust/rao:book} and the
references therein for a discussion on how the numerical radius of a
matrix can be estimated/computed.

Recall \cite{agler:ieot} that $T\in \LL(H)$ is called a $2$-hypercontraction
if
$$I - T^*T \geq 0 \quad \mbox{ and } \quad I - 2T^*T + T^{*2}T^2 \geq
0.$$

\begin{cor}
Let $T\in \LL(H)$ be a nilpotent $2$-hypercontraction with $T^n = 0$, $n
\geq 2$.
Then
$$\om_{\rho}(p(T)) \leq \om_{\rho}(p(B_{n}^*))$$
for all $\rho > 0$ and all polynomials $p$. Here $B_{n}^*\in \LL(\C^n)$ is
given by the matrix
$$B_{n}^* = \left(\begin{array}{cccccc}
            0 & \sqrt{\frac{1}{2}} & 0 & \cd & 0 & 0 \\
            0 & 0 & \sqrt{\frac{2}{3}} & \cd & 0 & 0 \\
            \vdots & \vdots & \vdots & \cd & \vdots & \vdots \\
            0 & 0 & 0 & \cd & 0 & \sqrt{\frac{n}{n+1}} \\
            0 & 0 & 0 & \cd & 0 & 0
        \end{array}\right). $$
\end{cor}

\begin{proof} Let $r < 1$ be a fixed positive real number.
Consider $\H = L^2_{a}(\D)$ the Bergman space of all analytic functions on $\D$
satisfying
$$\|f\|^2 = \frac{1}{\pi}\int_{\D}|f(re^{it})|^2\, dA < \infty ,$$
where $dA$ is the area Lebesgue measure. In this case $C(w,z) = (1 -
2wz + w^2z^2)^{-1}$ and $\H$ is a $\D$-space satisfying axiom (AEM)
with $\psi_{j}(z) = \sqrt{j+1}z^j$ and $f_{m}(w,z) = (m+1)w^mz^m -
mw^{m+1}z^{m+1}$. Then $Z^*$ is unitarily equivalent to the Bergman
shift $B^*$, where $B$ is given by
$Be_{p} = \sqrt{\frac{p+1}{p+2}}e_{p+1}$ for a
suitable orthonormal basis $(e_{p})$.

We have $\|rT\| \leq 1$, $(rT)^m \to 0$ strongly and also
\cite{agler:ieot}
$$I - 2(rT)^*(rT) + (rT)^{*2}(rT)^2 \geq 0.$$
It follows from \cite[Example 2]{aem} that $rT$ satisfies condition
(ii) of Theorem \ref{thm:unco}. It follows from Theorem \ref{thm:24}
that
$$\om_{\rho}(p(rT)) \leq \om_{\rho}(p(B_{n}^*)),$$
since $B_{n}^{*}$ is unitarily equivalent to $B^*\mid \Ker (B^{*n})$.
This holds for all $r < 1$ ; it also holds for $r = 1$.
\end{proof}

The numerical radius of $B_{n}^*$ can be expressed \cite{stout:pams}
in terms of the smallest positive root of a polynomial involving
circularly symmetric functions. To give a flavor of what can be done,
we prove here the
following inequalities.
\begin{cor}
Suppose $T \in \LL(H)$ satisfies $\|T\| \leq 1$, $T^3 = 0$ and
$I - 2T^*T + T^{*2}T^2 \geq  0$.
Then
$$\om_{2}(T) \leq \sqrt{\frac{7}{24}} \quad \mbox{ and } \quad
\om_{2}(T^2) \leq \sqrt{\frac{1}{12}}$$
and these constants are the best possible ones.
\end{cor}

\begin{proof} We have to compute $\om_{2}(B_{n}^*)$ and
$\om_{2}(B_{n}^{*2})$ for $n = 3$.
This can be
done using \cite{stout:pams} or in the following (equivalent) way.
Consider the symmetric $n\times n$ matrix
$$A_{n} = B_{n}^* + B_{n} = \left(\begin{array}{cccccc}
            0 & \sqrt{\frac{1}{2}} & 0 & \cd & 0 & 0 \\
            \sqrt{\frac{1}{2}} & 0 & \sqrt{\frac{2}{3}} & \cd &
0 & 0 \\
            \vdots & \vdots & \vdots & \cd & \vdots & \vdots \\
            0 & 0 & 0 & \cd & 0 & \sqrt{\frac{n}{n+1}} \\
            0 & 0 & 0 & \cd & \sqrt{\frac{n}{n+1}} & 0
        \end{array}\right). $$
Let $\theta$ be a real number. If $D(\theta)$ denotes the diagonal matrix with
$e^{ij\theta}$, $j = 1, \ldots , n$, on the main diagonal, then we
have
$D(\theta)^*(e^{i\theta}B_{n} + e^{-i\theta}B_{n}^*)D(\theta) = A_{n}$.
Recall that
$$\om_{2}(T) = \frac{1}{2}\sup_{\theta\in
\R}\|e^{i\theta}T +
e^{-i\theta}T^*\|.$$
Therefore
\begin{eqnarray*}
    \om_{2}(B_{n}^*) & = & \frac{1}{2}\sup_{\theta\in
\R}\|e^{i\theta}B_{n} +
e^{-i\theta}B_{n}^*\|  \\
     & = & \frac{1}{2}\sup_{\theta\in \R}\|D(\theta)^*(e^{i\theta}B_{n} +
     e^{-i\theta}B_{n}^*)D(\theta)\|  \\
     & = & \frac{1}{2}\|A_{n}\|.
\end{eqnarray*}
Since $\frac{1}{2}A_{n}$ is hermitian, its norm coincides with its
largest eigenvalue. For $n = 3$ it is equal to $\sqrt{7/24}$.
In a similar way, the numerical radius of $B_{3}^{*2}$ is the spectral
radius of $\frac{1}{2}(B_{3}^{*2} + B_{3}^{2})$, that is
$\sqrt{1/12}$. 
\end{proof}

Note that the inequality
$$\om_{2}(T^2) \leq \sqrt{\frac{1}{12}} = 0.2886\ldots$$
is
an improvement of the inequality
$$\om_{2}(T^2) \leq \om_{2}(T)^2 \leq
\frac{7}{24} = 0.2916\ldots \; .$$

\noi 

{\bf Inequalities for $n$-tuples of operators.} Theorem \ref{thm:24}
can be applied also for $n$-tuples of
commuting operators $T = (T_{1}, \cdots , T_{n}) \in
\LL(H)^n$, $n\geq 1$. In fact, anytime we dispose of a model
operator, the techniques of Section 2 can be used to obtain
constrained von Neumann inequalities.
We give only one example using the model of Vasilescu
\cite{vasilescu}. It corresponds, using the notation of Section 2,
to the domain
$$D = \{z\in \C^n : \sum_{j}c_{ij}|z_{j}|^2 < 1, 1\leq i \leq m\}.$$
This generalizes
previous models for the unit ball in $\C^n$
and for the unit polydisc (cf~. the references in \cite{vasilescu}).

Let $m \geq 1$ be a fixed integer. Let $p = (p_{1},
\cdots , p_{m})$ be a family of complex polynomials
$$p_{j}(z) = 1 - c_{j1}z_{1} - \ldots - c_{jn}z_{n}\; ,$$
for $j = 1, \cdots ,m$, $z = (z_{1}, \cdots , z_{n}) \in \C^n$ such
that
 \begin{itemize}
    \item $c_{jk} \geq 0$ for all indices $j$ and $k$ ;

    \item for every $k \in \{1, \cdots , n\}$ there is $j \in \{1,
\cdots ,
    m\}$ such that $c_{jk} \neq 0$

    \item $p_{j}$ is identical $1$ for no indice $j$.
 \end{itemize}
The case
$$p_{j}(z) = 1 - z_{j} \quad 1 \leq j \leq n\; ,$$
corresponds to the unit polydisc in $\C^n$, while
$$p_{1}(z) = 1 - z_{1} - \cdots - z_{n} $$
corrsponds to the unit ball.

If $\gamma = (\gm_{1}, \cdots , \gm_{m}) \in \Z^m_{+}$, we set
$$p^{\gm}(z) = p_{1}(z)^{\gm_{1}} \cdots p_{m}(z)^{\gm_{m}}
\quad (z \in \C^n)\; .$$
Define
$$V_{T,j} = \sum_{k=1}^n c_{jk}M_{T_{k}}, \quad j = 1, \cdots , m ,$$
and $V_{T} = (V_{T_{1}}, \cdots , V_{T_{n}})$. Define
$$\Delta_{T}^{\gm} = (I - V_{T,1})^{\gm_{1}} \cdots (I -
V_{T,m})^{\gm_{m}}(I_{H})\; ,$$
where $I_{H}$ is the identity on $H$ and $I = I_{\LL(H)}$
is the identity on $\LL(H)$.
Let $\gm \geq (1, \ldots , 1)$.
We say \cite{vasilescu} that $T\in \LL(H)^n$ satisfies the positivity
condition
$(p,\gamma)$ if
$$\Delta_{T}^{\beta} \geq 0 , \quad \mbox{ for all } \beta , 0 \leq
    \beta \leq \gm .$$
We denote by $S^{(p,\gamma)} \in \LL(\ell_{2}(\Z_{+}^n,\C))$
the \emph{backwards multishift} of type $(p,\gamma)$ as
defined in \cite{vasilescu} (in fact, $S^{(p,\gamma)}\otimes I_{H}$
is the model there).

\begin{cor}\label{cor:46}
    Suppose, with notation as above, that
    $\gm \geq (1, \ldots , 1)$.
    Let $T \in \LL(H)^n$ be a $n$-tuple of commuting operators satisfying
    the positivity condition $(p,\gm)$ and the constraint
    $q(T) = 0$ for a fixed polynomial $q$ in $n$ variables.
    Then, for any admissible
    radius $\om$ and any polynomial $f$ in $n$ variables, we have
    $$\om (f(T)) \leq \om (f(S^{(p,\gm)})\mid \Ker q(S^{(p,\gm)})))\; .$$
\end{cor}

\begin{proof} Let $r \in ]0,1[$. It was proved in
\cite[Proposition 3.15]{vasilescu} that $rT$ is unitarily equivalent
to the restriction of $S^{(p,\gamma)}\otimes I_{H}$ to an invariant
subspace. Using the admissibility of $\om$, the fact that $q(T) =
0$, and making $r\to 1$ at the end, we obtain the desired inequality.
\end{proof}

A proof of the above corollary can be given using directly Theorem
\ref{thm:24} (cf~. Example (5) in \cite{aem}).
The unconstrained von Neumann inequality in this case, for the operator
norm, is
\cite[Proposition 3.15]{vasilescu}.

\noi 

{\bf Applications of Theorem \ref{thm:25}.} The
following result is obtained from Theorem \ref{thm:25} in the
classical case $\H = H^2$.

\begin{cor}
Let $T\in \LL(H)$ be a Hilbert space contraction such that $Q(T^*,T) =
0$ for a given $Q \in P(\C)$ of degree $d$
with $Q(e^{-it},e^{it}) \neq 0$ for some $t \in
\R$. Then there exists an invariant subspace $E$ for the backward shift
$S^*$ on $H^2$ such that
$$\om_{\rho}(P(T^*,T)) \leq \om_{\rho}(P(S_{E},S^*_{E}))$$
for all $\rho \in ]0,2]$ and all $P\in P(\C)$. Here $S_{E}\in \LL(E)$ is 
the adjoint of $S_{E}^* = S^*\mid E$.
\end{cor}
It follows from the proof of Theorem \ref{thm:25} that the space $E$
in the above corollary is given by $E = \Ker S^{*d}Q(S,S^*)$.
The following is a possible application.
\begin{cor}
Let $m \geq n \geq 1$ be two positive integers.
Let $T\in \LL(H)$ be a contraction and suppose
that $T^{*m} = T^n$. Let $\rho\in ]0,2]$ and let $P \in P(\C)$.
Then
$$\om_{\rho}(P(T^*,T)) \leq \om_{\rho}(P(S_{m+n},S_{m+n}^*)) ;$$
in particular,
$$\om_{2}(T^l)) \leq \cos \frac{\pi}{[(m+n-1)/l]+2}$$
for all $l$ with $1\leq l \leq m+n-1$.
\end{cor}

\begin{proof} Set $Q(w,z) = w^m - z^n$. We have
$Q(e^{-it},e^{it}) \neq 0$ for some $t \in
\R$.
Note that $S$ on $H^2$ is unitarily
equivalent to the forward
shift on $S$ on $\ell_{2}$. We have
$S^{*m}Q(S,S^*) = I - S^{*(m+n)}$.
Then $E = \Ker S^{*m}Q(S,S^*)$ is given by
$$E = \{(h_{0}, h_{1}, \ldots, h_{p-1}, h_{0}, h_{1}, \ldots, 
h_{p-1}, h_{0}, \ldots) : h_{k}
\in \C \mbox{ for } 0\leq k \leq p\},$$
where $p = m+n$.
Thus
$S_{E}^* = S^*\mid E$ is unitarily equivalent to $S_{p}^{*(\infty)}$.
Since $\om_{\rho}$ is strongly admissible for $\rho \leq 2$, we obtain
$$\om_{\rho}(P(T^*,T)) \leq \om_{\rho}(P(S_{m+n},S_{m+n}^*)).$$
The second inequality is obtained for $P(w,z) = z^l$.
\end{proof}

\noi 

{\bf Applications of Theorem \ref{thm:62}.} Theorem \ref{thm:62}
can be applied for instance to bundles of the following
type.
Let $(p_{n})_{n\geq 0}$ be a sequence of polynomials which is uniformly
bounded on the closed unit disc. Suppose $v : \LL(H)\to [0,\infty[$ is such
that
$v(T)\left\| T\right\| \leq 1$.
Let $\mathcal{U}$ be a non trivial ultrafilter on $\N$.
For any $\alpha \in \D$ and any $T\in \LL(H)$, set
$$
R_{\alpha }(T)=\lim_{\mathcal{U}}[p_{n}(v(T)T)+p_{n}(v(T)T)^{\ast }].
$$
Denote $\mathcal{K} = \mathcal{K}_{R}$ the collection of operators
associated to the previous defined bundle $R$. Let $\om_{\mathcal{K}}$
be the associated operator radius.
With these notations and using
Theorem \ref{thm:62}, we obtain the following result.

\begin{cor}
    Assume that $u$ is a finite Blaschke product. Suppose that
    $T$ is a $C_{0}$ contraction such that $u(T)=0$ and $v(T) = v(S_{u})$.
    Then, with notation as above,
    $$
    \om_{\mathcal{K}}(f(T))\leq \om_{\mathcal{K}}(f(S(u)))
    $$
     for each $f\in A(\D)$.
\end{cor}

\section{Bounds of coefficients of positive rational
functions}\label{section:5}
There are many classical inequalities for coefficients of
(positive) trigonometric
polynomials. The next result shows the links between the numerical radius
of the extremal operator in the constrained von Neumann
inequalities and the Taylor coefficients of rational functions positive
on $\T$.

\begin{theorem}\label{thm:rat}
 Let $F=P/Q$ be a rational function
with no principal part and which is positive on the torus. Then the Taylor
coefficient $c_{k}$ of order $k$ satisfies the following
inequality
$$
|c_{k}| \leq c_{0}\om_{2}(R^{k}),
$$
where $R=S^{\ast }\mid \Ker (Q(S^{\ast }))$.
\end{theorem}

\begin{proof} First, observe that by continuity we may
assume that $F$ is strictly positive on the torus. Let $F=P/Q$ be a rational
function without principal part, that is we have
$d{{}^\circ}(P)<d{{}^\circ}(Q)$ for the degrees.
Assume that $F(z)>0$ for every $z\in \T .$ Denote by $\beta
_{1},...,\beta _{q}$ the zeros of $Q$ which are contained in the open unit
disc $\D$ and write $Q(z)=(z-\beta _{1})^{d_{1}}...(z-\beta
_{q})^{d_{q}}Q_{2}(z)$, where $Q_{2}$ has no zero in $\D$.
Consider the function $G(z)=\overline{F(1/\overline{z})}$ which is
analytic,
except at a finite set of complex numbers. Since $F$ is real on the torus,
we have $G(e^{it})=\overline{F(e^{it})}=F(e^{it})$ for every $t\in \R$.
The analytic extension principle implies that $F(z)=G(z)$ except for a
finite set in $\C$. Thus $F(z)$ can be written in the
following way
\begin{equation*}
F(z)=\frac{P(z)}{Q_{1}(z)Q_{2}(z)},
\end{equation*}
where $Q_{1}(z)=(z-\beta _{1})^{d_{1}}...(z-\beta _{q})^{d_{q}}$ and $
Q_{2}(z)=(1-\overline{\beta _{1}}z)^{d_{1}}...(1-\overline{\beta _{q}}
z)^{d_{q}}$. Because of the condition $F(z)=\overline{F(1/\overline{z
})}$, we have $P(z)=z^{2d}\overline{P(1/\overline{z})}$ where $
d=d_{1}+...+d_{q}=d
{{}^\circ}(Q)/2$. If $P(\alpha )=0$, with $\alpha \neq 0$,
then necessarily
$P(1/\overline{\alpha })=0$. Therefore $P$ can be written as
\begin{equation*}
P(z)=cz^{m_{0}}(z-\alpha _{1})^{m_{1}}...(z-\alpha _{p})^{m_{p}}(1-\overline{
\alpha _{1}}z)^{m_{1}}...(1-\overline{\alpha _{p}}z)^{m_{p}}
\end{equation*}
with a suitable constant $c$. We have $d=m_{1}+...+m_{p}$. Finally, we get
\begin{equation*}
F(e^{it})=c\left| \frac{P_{1}(e^{it})}{Q_{2}(e^{it})}\right| ^{2}
\end{equation*}
with $d{{}^\circ}
(P_{1})<d
{{}^\circ}
(Q_{2})$ and $c>0$. Note that
\begin{equation*}
\frac{P_{1}(z)}{Q_{1}(z)}=\sum_{k=1}^{q}\sum_{i=1}^{m_{k}}\frac{a_{k,i}}{(1-
\overline{\alpha _{k}}z)^{i}}
\end{equation*}
for some $a_{k,i}\in \C$. It follows that $P_{1}(z)/Q_{1}(z)\in
E:=H^{2}\ominus bH^{2}$, where $b$ is the associated Blaschke product defined
by
\begin{equation*}
b(z)=\prod_{k=1}^{q}(\frac{z-\alpha _{k}}{1-\overline{\alpha _{k}}z}
)^{m_{k}}.
\end{equation*}
It follows from Lemma \ref{lem:63}
that we have $F=\left| f\right| ^{2}$ with a suitable $f\in E$.

Denote by $R$ the restriction of the backward shift $S^{\ast }$ to the
invariant subspace $\Ker Q (S)^{\ast }$. Then, for any integer $k$, we get
\begin{equation*}
|c_{k}| =\left| <R^{k}f\mid f>\right| \leq \om_{2}(R^{k})\left\|
f\right\| _{2}^{2}=\om_{2}(R^{k})\left\| F\right\| _{1}=\om_{2}(R^{k})c_{0}.
\end{equation*}
This ends the proof. 
\end{proof}

Setting $Q(z)=z^{n-1}$ in the previous theorem, and using previous
computations of the numerical radii, we obtain the following
classical inequality due to E.V.~Egerv\'ary and O.~Sz\'azs (1927).
The bound for
$c_{1}$ is due to L. Fejer (1915).
\begin{cor}[Egerv\'ary-Sz\'azs]
Let $P(e^{it})=\sum_{j=-n+1}^{n-1}c_{j}e^{ijt}$
be a positive trigonometric polynomial ($n\geq 2$). Then
\begin{equation*}
\left| c_{k}\right| \leq c_{0}\cos (\frac{\pi }{[\frac{n-1}{k}]+2})\quad
\text{for }1\leq k\leq n-1.
\end{equation*}
\end{cor}
\begin{remark}
We note the amuzing consequence that Fejer's inequality for $\left|
c_{1}\right|$ implies, via operator inequalities, the Egerv\'ary-Sz\'azs
inequality. Indeed, by \cite{HaHa}, Fejer's inequality implies the
Haagerup-de la Harpe inequality (\ref{HHineq}). By Remark \ref{rem:42}
this implies a bound for $\om_{2}(T^m)$, which in turn implies, as in
\cite{HaHa}, the Egerv\'ary-Sz\'azs
inequality.
\end{remark}

The next result gives estimates involving two coefficients of a
positive trigonometric polynomial.

\begin{theorem}\label{thm:harm} Let $P(e^{it})=\sum_{j=-n+1}^{n-1}c_{j}e^{ijt}$
be a positive trigonometric polynomial ($n\geq2$). Then,
for every distinct numbers $k$ and $l$ among $\{0;...,n-1\}$, there
exists $\gamma \in \R$ such that
\begin{equation*}
\left| c_{k}\right| +\left| c_{l}\right| \leq c_{0}
\om_{2}(S_{n}^{k}+e^{i\gamma }S_{n}^{l}).
\end{equation*}
In particular, we have
\begin{equation*}
\left| c_{k}\right| +\left| c_{l}\right| \leq c_{0}\left(1+\cos \frac{\pi }{[
\frac{n-1}{k+l}]+2}\right)^{1/2}
\left(1+\cos \frac{\pi }{[\frac{n-1}{\left| k-l\right|
}]+2}\right)^{1/2}.
\end{equation*}
\end{theorem}

\begin{proof} We can assume that $c_{0}=1$. Since $P$ is positive,
we have $P=\left| Q\right| ^{2}$ for some $Q\in
\C_{n-1}[X]$, the space of all
polynomials of degree less or equal to $n-1$. For any $
k,l$, there exists $\gamma $ such that
\begin{equation*}
\left| c_{k}\right| +\left| c_{l}\right| =\left| c_{k}+e^{i\gamma
}c_{l}\right| =\left| \int_{0}^{2\pi }(e^{ik\theta }+e^{i(l\theta +\gamma
)})\left| Q(e^{il\theta })\right| ^{2}dm(\theta )\right| .
\end{equation*}
We deduce from the equality $\left\| Q\right\| _{2}=c_{0}=1$
that
\begin{equation*}
\left| c_{k}\right| +\left| c_{l}\right| \leq \om_{2}(S_{n}^{k}+e^{i\gamma
}S_{n}^{l}).
\end{equation*}
Denote $M = \om_{2}(S_{n}^{k}+e^{i\gamma
}S_{n}^{l})$. We have
\begin{eqnarray*}
M &=&\sup_{\left\| R\right\| _{2}\leq
1}\sup_{\alpha \in \R}\mbox{Re }e^{i\omega }\int_{0}^{2\pi
}(e^{ik\theta }+e^{i(l\theta +\gamma )})\left| R(e^{il\theta })\right|
^{2}dm(\theta ) \\
&=&2\sup_{\left\| R\right\|_{2}\leq 1}\sup_{\alpha \in \R%
}\int_{0}^{2\pi }\cos (\frac{1}{2}[(k+l)\theta +\gamma +2\alpha ])\cos (
\frac{1}{2}[(k-l)\theta -\gamma )]) \\
& & \times\left| R(e^{il\theta })\right|
^{2}dm(\theta ) \\
 &\leq &2\left(\sup_{\left\| R\right\| _{2}\leq 1}\sup_{\alpha \in \R%
}\int_{0}^{2\pi }\cos ^{2}(\frac{1}{2}[(k+l)\theta +\gamma +2\alpha ])\left|
R(e^{il\theta })\right| ^{2}dm(\theta )\right)^{1/2}\times  \\
&&\left(\sup_{\left\| R\right\| _{2}\leq 1}\int_{0}^{2\pi }\cos ^{2}(\frac{1%
}{2}[(k-l)\theta -\gamma ])\left| R(e^{il\theta })\right| ^{2}dm(\theta
)\right)^{1/2}.
\end{eqnarray*}
Let $R$ be in $\C_{n-1}[X]$ with $\left\| R\right\| _{2}\leq 1$. Since $%
L(e^{it})=R(e^{i(t-\frac{\gamma +2\omega }{k+l})})$ is also in
$\C_{n-1}[X]$ and of norm less or equal to one,
we obtain, using the rotation invariance of the Haar
measure, that
\begin{eqnarray*}
&&\sup_{\alpha \in \R}\sup_{\left\| R\right\| _{2}\leq
1}\int_{0}^{2\pi }\cos ^{2}(\frac{1}{2}[(k+l)\theta +\gamma +2\alpha %
])\left| R(e^{il\theta })\right| ^{2}dm(\theta ) \\
&=&\sup_{\left\| L\right\| _{2}\leq 1}\int_{0}^{2\pi }\cos ^{2}((\frac{k+l}{2%
})t)\left| L(e^{ilt})\right| ^{2}dm(\theta ) \\
&=&\frac{1}{2}+\frac{1}{2}\sup_{\left\| L\right\| _{2}\leq 1}\int_{0}^{2\pi
}\cos ((k+l)t)\left| L(e^{ilt})\right| ^{2}dm(\theta )=\frac{1}{2}%
(1+\om_{2}(S_{n}^{k+l})).
\end{eqnarray*}
In a similar way
\begin{equation*}
\sup_{\left\| R\right\| _{2}\leq 1}\int_{0}^{2\pi }\cos ^{2}(\frac{1}{2}[%
(k+l)\theta -\gamma ])\left| R(e^{il\theta })\right| ^{2}dm(\theta )=\frac{1%
}{2}(1+\om_{2}(S_{n}^{k-l})).
\end{equation*}
Finally, we obtain
\begin{equation*}
\left| c_{k}\right| +\left| c_{l}\right| \leq \sqrt{1+\om_{2}(S_{n}^{k+l})}%
\sqrt{1+\om_{2}(S_{n}^{k-l})}.
\end{equation*}
Since $\om_{2}(S_{n}^{p})=\cos (\frac{\pi }{[\frac{n-1}{p}]+2})$, we get the
desired result.
\end{proof}
\begin{remark}
\firsta For $l = 0$ we reobtain the Egerv\'ary-Sz\'azs
inequality.

\condb When $k+l > n-1$, we get from Theorem \ref{thm:harm} that
\begin{equation*}
\left| c_{k}\right| +\left| c_{l}\right| \leq c_{0}\left(1+\cos \frac{\pi }{[%
\frac{n-1}{\left| k-l\right| }]+2}\right)^{1/2}.
\end{equation*}
In particular, if $n\geq 4$, we obtain
\begin{equation*}
\left| c_{1}\right| +\left| c_{n-1}\right| \leq c_{0}\sqrt{3/2}.
\end{equation*}
This estimate is better than that one obtained by applying twice the
Egerv\'ary-Sz\'azs inequality.

\condc In some particular cases, it is possible to compute exactly the
numerical radius $M = \om_{2}(S_{n}^{k}+e^{i\gamma
}S_{n}^{l})$. Suppose $n = 9$, $k = 3$, $l = 7$.
It follows from \cite{davidson/holbrook} that $M =
\cos (\pi/10)$ if $\gamma =
0$. The method from \cite{davidson/holbrook}
does not seem to apply for an arbitrary $\gamma$.
\end{remark}

\section{Other type of constraints} The constraints until
now were of algebraic type ($q(T) = 0$ or
$Q(T^*,T) = 0$). We discuss briefly constraints of different nature.

\noi 

{\bf Some positivity conditions.} We discuss constrained von
Neumann inequalities for the numerical radius $\om_{2}$ of an operator
satisfying some positivity conditions $R_{\lambda}(T^*,T) \geq 0$ for
$\lambda \in \T$.

\begin{prop}\label{es}
Let $n \geq 2$ be a positive integer and let $\rho_{k}$, $0\leq k \leq
n-1$, be $n$ positive reals with $\rho_{0} = 1$. Let $T \in \LL(H)$ be
an operator such
that $R_{\lambda}(T^*,T) \geq 0$ for
$\lambda \in \T$, where
\begin{equation}
    R_{\lambda}(w,z) = 1 +
    \sum_{k=1}^{n-1}\frac{\lambda^k}{\rho_{k}}w^k +
    \sum_{k=1}^{n-1}\frac{\overline{\lambda}^k}{\rho_{k}}z^k \;
    (\lambda \in \T).
    \label{rho/sections}
\end{equation}
Then
$$\om_{2}(T^m) \leq \rho_{m}\cos \frac{\pi}{[\frac{n-1}{m}]+2}$$
for each $m \in \{1, 2, \ldots, n-1\}$.
\end{prop}

\begin{proof} Let $h \in H$ be a norm-one
vector and let $\theta \in \R$. Set
$$c_{k} \;  = \left\{
\begin{array} {r@{\quad:\quad}l} 1 & \mbox{if } k = 0 \\
\frac{1}{\rho_{k}}\langle T^kh|h \rangle & \mbox{if } k > 0 \\
\frac{1}{\rho_{|k|}}
\langle h|T^{|k|}h \rangle & \mbox{if } k < 0 \end{array} \right. $$
and
$$t_{n}(\theta) = \sum_{k=-n+1}^{n-1}c_{k}e^{ik\theta} \; .$$
Then $t_{n}$ is a positive trigonometric polynomial since
$$t_{n}(\theta) = \langle R_{\exp(it)}(T^*,T)h|h \rangle .$$
According to the Egerv\'ary-Sz\'azs inequality, we have
$$\frac{1}{\rho_{m}}|\langle T^mh|h \rangle | = |c_{m}| \leq \cos (\frac{\pi
}{[\frac{n-1}{m}]+2})$$
which gives the desired inequality. 
\end{proof}

If $\rho_{k} = 1$, for each $k\leq n-1$, then $R_{\lambda}(T^*,T)$ in
Equation (\ref{rho/sections}) are the $n$th sections of the operator
kernel $K_{\lambda}(T)$.

In fact, the following more general result holds.
\begin{theorem}\label{thm:662}
Let $n \geq 2$ be a positive integer and let $\rho_{k}$, $0\leq k \leq
n-1$, be $n$ positive reals with $\rho_{0} = 1$. Let $T \in \LL(H)$ be
an operator such
that $R_{\lambda}(T^*,T) \geq 0$ for
$\lambda \in \T$, where $R_{\lambda}(w,z)$ are given by
(\ref{rho/sections}). Then, for any strongly admissible radius $\nu$
and any $m \in \{0, 2, \ldots, n-1\}$,
we have
$$\nu (T^m) \leq \rho_{m}\nu (S_{n}^{*m}).$$
Moreover, if $m$ and $l$ are distinct numbers among $\{0;...,n-1\}$
and if additionally $\rho_{m} = \rho_{l}$, then
$$\nu (T^m+T^l) \leq \rho_{m}\nu (S_{n}^{*m}+S_{n}^{*l}).$$
In particular, we have
$$
\omega _{\rho }(T^{m}+T^{l})\leq 
\left(\frac{2}{\rho}-1\right)2^{\rho }\rho _{m}\left[
1+\cos \frac{\pi}{[\frac{n-1}{m+l}]+2}\right]^{\frac{1-\rho}{2}}\left[
1+\cos \frac{\pi}{[\frac{n-1}{\left|m-l\right|}]+2}\right]^{\frac{1-\rho
}{2}}$$ 
if $\rho \in ]0,1]$, and 
$$
\omega _{\rho }(T^{m}+T^{l})\leq 2^{2-\rho }\rho _{m}^{\rho -1}\left[ 1+\cos
\frac{\pi }{[\frac{n-1}{m+l}]+2}\right] ^{\frac{\rho -1}{2}}\left[ 1+\cos
\frac{\pi }{[\frac{n-1}{\left| m-l\right| }]+2}\right] ^{\frac{\rho -1}{2}}
$$
if $\rho \in ]1,2]$.
\end{theorem}

Its proof follows from Theorem \ref{thm:harm},
interpolation properties of $\omega _{\rho }$ 
(see \cite[p. 296]{fong/holbrook}) and
the following generalization of a result of W. Arveson
(obtained in \cite{arveson} for $\rho_{k} = 1$, $k\geq 1$).

\begin{theorem} \label{arv}
    Let $T \in \LL(H)$ be a contraction and let $n \geq 2$.
    Suppose $T$ satisfies $R_{\lambda}(T^*,T) \geq 0$ for all $\lambda
        \in \T$, where $R_{\lambda}(w,z)$ are given by
(\ref{rho/sections}).
Then there is a Hilbert space $K \supset H$ and a nilpotent
        contraction
        $N \in \LL(K)$ such that $N^n = 0$, $N$ is unitarily
        equivalent to $S_{n}^{*(d)}$, $d$ finite or $\infty$, and
        $T^j = \rho_{j}P_{H}N^j\mid
        H$ for $j = 0, 1, \cdots , n-1$.
\end{theorem}

\begin{proof} The idea of the proof is that of \cite{arveson} and
some details will be omitted below. Define a
linear map $\varphi$ from $span \{S_{n}^{*j} : 0 \leq j \leq n-1\}$
 onto $span \{\frac{1}{\rho_{j}}T^j :
0 \leq j \leq n-1\}$ by
$\varphi (S_{n}^{*j}) = \frac{1}{\rho_{j}}T^j$ and by linearity.
Define the map $\psi : C(\T) \to \LL(H)$ by
$$\psi(f) = \frac{1}{2\pi}\int_{0}^{2\pi}f(e^{i\theta})\,
R_{\exp(i\theta)}(T^*,T)\, d\theta \ .$$
It is a positive linear map.
Note that $\psi(z^j) = \frac{1}{\rho_{j}}T^j$ for $j = 0, 1, \cdots , n-1$ and
$\psi(z^j) = 0$ for $j \geq n$, where $z(\theta) = \theta$. It is
known that a positive
linear map on a commutative C${}^{\ast}$-algebra is completely
positive and a completely positive map which preserves the identity is
completely contractive \cite{paulsen:book}.
The restriction $\psi_{0}$ of $\psi$ on the
disc algebra (the closed linear span in $C(\T)$ of $1,z, z^2,
\cdots $) is a completely contractive linear map such that
$\psi_{0}(z^j) = \frac{1}{\rho_{j}}T^j$ for $j = 0, 1, \cdots , n-1$ and
$\psi_{0}(z^j) = 0$ for $j \geq n$. It vanishes on the ideal $z^nA$
and thus it induces a completely contractive linear map
$\psi_{00}$ of the quotient $A/z^nA$ into $\LL(H)$. It
was proved in \cite{arveson}
that $\mu (S_{n}^{*j}) = z^j + z^nA$ defines a completely
isometric linear map of $span \{I, S_{n}^*, \cdots , S_{n}^{*n-1}\}$
onto $A/z^nA$. The original map $\varphi = \psi_{00}\mu$ is
thus completely contractive. Since $\varphi (I) = I$, $\varphi$ has
\cite{paulsen:book} a
completely positive extension to $C^{\ast}(S_{n}^*) =
\LL(\C^n)$. Stinespring's theorem \cite{paulsen:book} furnishes then
a unital C${}^*$-representation $\pi$. Then $N = \pi(S_{n}^*)$
gives, as in \cite{arveson}, the desired representation. 
\end{proof}

In the case $\rho_{k} = 1$ for all $k$, studied in \cite{arveson}, 
the converse of Theorem \ref{arv} also holds.
Also, an operator $T$ satisfies
$$I + 2 \mbox{ Re }\sum_{k=1}^{n-1}z^kT^k \geq 0, \mbox{ for each }
  z \in \T ,$$
if and only if \cite{arveson}
$$2 \mathrm{Re }(I - zT)^{\ast}z^nT^n \leq I - T^{\ast}T , \mbox{ for each }
  z \in \T \; .$$
In particular this holds if $T$ satisfies $T^n = 0$ and $I -
T^{\ast}T \geq 0$.

\noi 

{\bf Stability of the algebraic constraints.} In what follows 
$\ep > 0$ is supposed to be a (fixed) small positive number.
We study what
happens if the constraint $q(T) = 0$ is
replaced by $\|q(T)\| \leq \ep$.

\begin{prop}
Let $q$ be a polynomial. For each $\ep > 0$
there exists  $\delta > 0$ such that every contraction $T\in \LL(H)$
with $\|q(T)\| \leq \delta$ satisfies
$$\om_{\rho}(T) \leq \ep + \om_{\rho}(S^*|\Ker q(S^*))$$
for every $\rho\in ]0,2]$.
\end{prop}

\begin{proof} By \cite[Corollary 2.22]{herreroI:book}, for
every $\ep > 0$ there is $\delta > 0$ such that, if $\|T\| \leq 1$
and $\|q(T)\| \leq \delta$, then there exists $T'\in \LL(H)$ such
that $q(T') = 0$ and $\|T - T'\| \leq \ep$. Note also that
$\om_{\rho}$ is a norm for $\rho \leq 2$. We thus have
\begin{eqnarray*}
    \om_{\rho}(T) & \leq & \om_{\rho}(T-T') + \om_{\rho}(T')   \\
     & \leq & \|T-T'\| + \om_{\rho}(T')\\
     & \leq & \ep + \om_{\rho}(S^*|\Ker q(S^*)).
\end{eqnarray*}
The proof is complete. 
\end{proof}

It was proved in \cite[Theorem 2.21]{herreroI:book} that if $\|T\| \leq 1$
and $\|T^n\| \leq \ep$, then there exists $T'\in \LL(H)$ such
that $T^{'n} = 0$ and $\|T - T'\| \leq \delta_{n}(\ep)$, where
$\delta_{n}(\ep)$ is defined inductively by
$$\delta_{2}(\ep) =
(2\ep)^{1/2} \quad \mbox{ and } \quad \delta_{k}(\ep) = \{\ep +
[\delta_{k-1}((k-1)\sqrt{\ep})]^2\}^{1/2}.$$
This implies that if $\|T\| \leq 1$
and $\|T^n\| \leq \ep$ then
$$\om_{2}(T) \leq \cos (\frac{\pi
}{n+1}) + \delta_{n}(\ep) .$$
Note that $\lim_{\ep\to 0}\delta_{n}(\ep) = 0$.

The following result gives a better bound for small $\ep$; we obtain the
Haagerup-de la Harpe inequality for $\ep \to 0$.

\begin{theorem}\label{epsilonized}
Let
$n \geq 2$ be a positive integer.
Suppose $T\in \LL(H)$ is a
contraction satisfying $\|T^n\| \leq \ep$ and
$$\sum_{k> n+1}\|T^k\| < + \infty .$$
Then
\begin{eqnarray*}
    \om_{2}(T) & \leq & \cos (\frac{\pi
}{n+1}) + 3\left[\pi\cos^4\frac{\pi}{2(n+1)}\right]^{1/3}\left(
\frac{\ep}{n+1}\right)^{2/3}  \\
     & \leq & \cos (\frac{\pi
}{n+1}) + 3\sqrt[3]{\pi}\left(
\frac{\ep}{n+1}\right)^{2/3}.
\end{eqnarray*}
\end{theorem}

The proof uses the following epsilonized Fejer
inequality. Note that an epsilonized version of the
Egerv\'ary-Sz\'azs inequality can be proved along the same lines.
\begin{lem}[The epsilonized Fejer inequality]
Let $h$ be a positive function,
$$h(\theta) = \sum_{m\in \Z}c_{m}e^{im\theta},$$
such that $\sum_{m\in\Z}|c_{m}| < \infty$
with $c_{0} = 1$ and $|c_{k}| \leq \ep$ for
$k \geq n$. Then
$$|c_{1}| \leq \cos (\frac{\pi
}{n+1}) + 3\left[\pi\cos^4\frac{\pi}{2(n+1)}\right]^{1/3}\left(
\frac{\ep}{n+1}\right)^{2/3}.$$
\end{lem}

\begin{proof} The following result has been proved in
\cite[Example 4(a)]{janssen}~:
Let $f$ be the Fourier transform of a non-negative integrable
function $\varphi$~:
$$f(x) = \int_{-\infty}^{\infty}e^{ixt}\varphi(t)\, dt .$$
Let $u > 0$
and suppose that $f(0) = 1$ and $|f(ku)| \leq \ep$ for $k\geq n$.
Then
$$|f(u)| \leq \cos (\frac{\pi
}{n+1}) + 3\left[\pi\cos^4\frac{\pi}{2(n+1)}\right]^{1/3}\left(
\frac{\ep}{n+1}\right)^{2/3} .$$
This is a generalization of a result due to Boas and Kac \cite{BoKa}
for band-limited functions.

Set now $\varphi (t) = h(-t)$, $t\in [-\pi,\pi]$. Consider $f$
the Fourier transform of $\varphi$. Then $f(0) = c_{0} = 1$, $f(k) = c_{k}$
and thus $|f(k)| \leq \ep$ for $k \geq n+1$.
We can now apply \cite{janssen} with $u = 1$. 
\end{proof}

\begin{proof}[Proof of Theorem \ref{epsilonized}] The proof is similar
to the proof of Proposition \ref{es}. By replacing eventually $T$ by
$rT$, $0 < r < 1$, it is possible to assume that the spectrum of $T$
is contained in $\D$. For each norm-one
vector $h \in H$ and $\theta \in \R$, set
$$c_{k} \; (= c_{k}(h)) = \left\{
\begin{array} {r@{\quad:\quad}l} 1 & \mbox{if } k = 0 \\
\langle T^kh|h \rangle & \mbox{if } k > 0 \\
\langle h|T^{|k|}h \rangle & \mbox{if } k < 0 \end{array} \right. $$
and
$$h(\theta) = \sum_{k\in \Z}c_{k}e^{ik\theta} \; .$$
Then $\sum_{m\in\Z}|c_{m}| < \infty$.
Note also that
$$h(\theta) = \langle K_{\exp(it)}(T)h|h \rangle $$
and the operator kernel
$$K_{\exp(it)}(T) = (I-e^{it}T^*)^{-1}(I - T^*T)(I-e^{-it}T)^{-1}$$
is positive since $T$  is a contraction. We use now the epsilonized
Fejer inequality. 
\end{proof}

\begin{cor}
Let $n$ and $m$ be two positive integers such that
$m\geq n \geq 2$.
Suppose $T\in \LL(H)$ is a
contraction satisfying $\|T^n\| \leq \ep$ and $T^m = 0$. Then
$$\om_{2}(T) \leq \min \left[\; \cos (\frac{\pi
}{m+1}) \; ; \;  \cos (\frac{\pi
}{n+1}) + 3\sqrt[3]{\pi}\left(
\frac{\ep}{n+1}\right)^{2/3}\; \right] .$$
\end{cor}

\medskip

\noi {\bf Acknowledgment}. 
This paper was written during several
visits of the authors to the Mathematical Departements of
Universities of Lille and Lyon. We wish to
thank both institutions for their help.

\end{document}